\documentclass[11pt]{amsart}
\usepackage[T1]{fontenc}
\usepackage[utf8]{inputenc}
\usepackage[british]{babel}
\usepackage[margin=1in]{geometry}
\usepackage{amsmath,amssymb,amsthm}
\usepackage{palatino,euler,microtype}

\usepackage{booktabs,longtable,array}
\usepackage{tikz}
\usetikzlibrary{positioning,fit,backgrounds,calc,shapes.geometric}
\DeclareMathVersion{jfefpalatino}
\SetSymbolFont{operators}{jfefpalatino}{OT1}{zplm}{m}{n}
\SetSymbolFont{letters}{jfefpalatino}{OML}{zplm}{m}{it}
\SetSymbolFont{symbols}{jfefpalatino}{OMS}{zplm}{m}{n}
\SetSymbolFont{largesymbols}{jfefpalatino}{OMX}{zplm}{m}{n}
\SetMathAlphabet{\mathbf}{jfefpalatino}{OT1}{zplm}{b}{n}
\SetMathAlphabet{\mathit}{jfefpalatino}{OT1}{ppl}{m}{it}
\tikzset{
  jfef text family/.initial=ppl,
  jfef math version/.initial=jfefpalatino,
  jfef font size/.initial=8,
  jfef baseline/.initial=9,
  jfef label font/.style={font={
    \fontfamily{\pgfkeysvalueof{/tikz/jfef text family}}\selectfont
    \fontsize{\pgfkeysvalueof{/tikz/jfef font size}}{\pgfkeysvalueof{/tikz/jfef baseline}}\selectfont
    \mathversion{\pgfkeysvalueof{/tikz/jfef math version}}
  }},
  jfef vertex/.style={
    circle,draw=black!70,line width=.4pt,
    minimum size=5.8mm,inner sep=.8pt,outer sep=0pt,
    fill=white,opacity=1,fill opacity=1,draw opacity=1,text opacity=1,
    text=black,jfef label font
  },
  jfef type/.style={jfef vertex,minimum size=6mm,inner sep=1pt},
  jfef centre/.style={jfef vertex,minimum size=8mm,inner sep=1pt,fill=orange!20},
  jfef auxiliary/.style={jfef vertex,rectangle,minimum size=8mm,inner sep=1pt,fill=cyan!15},
  jfef edge/.style={draw=black!27,line width=.25pt,line cap=round},
  jfef complement edge/.style={draw=black!48,line width=.35pt,line cap=round},
  jfef new edge/.style={draw=blue!55!black,opacity=.55,line width=.35pt,line cap=round},
  jfef new complement edge/.style={draw=blue!65!black,line width=.6pt,line cap=round}
}

\usepackage{float,placeins,needspace,enumitem,siunitx}
\usepackage{xcolor}
\definecolor{darkblue}{rgb}{0,0,0.7}
\definecolor{darkred}{rgb}{0.7,0,0}
\newcommand{\defi}[1]{\textit{\color{blue}#1}}
\usepackage[numbers]{natbib}
\usepackage{hyperref}
\hypersetup{colorlinks=true,linkcolor=darkred,citecolor=darkblue,urlcolor=darkblue,
  pdfauthor={Jesús F. Espinoza, Martín Eduardo Frías-Armenta, Rosalía G. Hernández-Amador, Héctor A. Hernández-Hernández},
  pdftitle={Graph reductions and the topology of flag complexes}}
\newcommand{\I}{\mathfrak I}
\newcommand{\IS}{\mathfrak I_S}
\newcommand{\ISV}{\mathfrak I_{SV}}
\newcommand{\Dis}[1]{#1\text{-Dis}}
\newcommand{\ZZ}{\mathbb Z}
\newcommand{\Ztwo}{\mathbb F_2}
\newcommand{\cl}{\operatorname{cl}}
\newcommand{\Nb}{N}
\newcommand{\Gt}{G_{13}}
\newcommand{\Ga}{G_{14}^{a}}
\newcommand{\Gb}{G_{14}^{b}}
\newcommand{\Gc}{G_{14}^{c}}
\newcommand{\Dunce}{D_{15}}
\newtheorem{theorem}{Theorem}[section]
\newtheorem{lemma}[theorem]{Lemma}
\newtheorem{proposition}[theorem]{Proposition}
\newtheorem{corollary}[theorem]{Corollary}
\newtheorem{definition}[theorem]{Definition}

\theoremstyle{definition}

\newcommand{\RowThirteenAcyclic}{\num{5619022190286}}

\newcommand{\RowThirteenIS}{\num{5619022190283}}
\newcommand{\RowFourteenAcyclic}{\num{755877673752}}

\newcommand{\RowFourteenIS}{\num{755877673119}}
\newcommand{\RowFourteenCone}{\num{20797002}}
\newcommand{\NumTwoDisNotOneDisThirteen}{\num{35757}}

\title[Graph reductions and the topology of flag complexes]{Graph reductions and the topology of flag complexes}
\author[J. F. Espinoza]{Jes\'us F. Espinoza\textsuperscript{*}}
\thanks{\textsuperscript{*}Corresponding author}
\author[M. E. Fr\'ias-Armenta]{Mart\'in Eduardo Fr\'ias-Armenta}
\author[R. G. Hern\'andez-Amador]{Rosal\'ia G. Hern\'andez-Amador}
\author[H. A. Hern\'andez-Hern\'andez]{H\'ector~A.~Hern\'andez-Hern\'andez}
\address{Departamento de Matem\'aticas, Universidad de Sonora,
  Edificio 3K-1, Blvd. Luis Encinas y Rosales s/n, Colonia Centro,
  Hermosillo, Sonora, 83000, M\'exico}
\email[Jes\'us F. Espinoza]{jesusfrancisco.espinoza@unison.mx}
\urladdr[Jes\'us F. Espinoza]{https://orcid.org/0000-0001-5601-3634}
\email[Mart\'in Eduardo Fr\'ias-Armenta]{eduardo.frias@unison.mx}
\urladdr[Mart\'in Eduardo Fr\'ias-Armenta]{https://orcid.org/0000-0002-7687-1270}
\email[Rosal\'ia G. Hern\'andez-Amador]{rosalia.hernandez@unison.mx}
\urladdr[Rosal\'ia G. Hern\'andez-Amador]{https://orcid.org/0000-0001-9657-7012}
\email[H\'ector A. Hern\'andez-Hern\'andez]{hector.hernandez@unison.mx}
\urladdr[H\'ector A. Hern\'andez-Hern\'andez]{https://orcid.org/0000-0002-4046-4841}
\date{}
\subjclass[2020]{Primary 05C75; Secondary 05C30, 05C85, 05E45, 55U10, 57Q10}
\keywords{contractible transformations, strong I-contractibility, dismantlability,
  clique complex, odd torsion, enumeration}
\begin{document}
\begin{abstract}
Ivashchenko's transformations characterise the graphs whose clique
complexes are contractible. A graph is strong I-contractible if it
reduces to one vertex by deleting vertices and edges. At each step the
neighbours of the vertex, or the common neighbours of the endpoints
of the edge, must induce a strong I-contractible graph. Write $\IS$ for this family. We show
that $14$ is the least order of an I-contractible graph outside $\IS$
and also the least number of vertices of a contractible, non-collapsible
flag complex. Our example attains both bounds. Its clique complex has
no free face, yet adding one edge to the graph gives a member of
$\IS$. Since the common neighbours of its endpoints induce a path,
the new edge can be deleted admissibly. If
contractibility is weakened to acyclicity over $\mathbb F_2$, the
minimum order outside $\IS$ drops to $13$. There are exactly three
examples of that order up to isomorphism, each with first integral
homology $\mathbb Z/3$. They also attain the least number of vertices
of a flag complex with odd torsion in first homology. These three
graphs disprove the acyclicity
assertion of Conjecture~14 of Dochtermann, Espinoza, Fr\'ias-Armenta and
Hern\'andez. The graph of order $14$ and its one-edge extension
disprove the remaining assertion and Conjecture~15, respectively. The
lower bounds follow from an enumeration of connected graphs through
order $13$. Among $K_4$-free graphs of order $14$, exactly eight
isomorphism classes are I-contractible and lie outside $\IS$.
\end{abstract}
\maketitle
\enlargethispage{2pt}
\section{Introduction}\label{sec:intro}

A finite \defi{flag complex} is determined by its graph of vertices and
edges, since its simplices are precisely the cliques of that graph.
This correspondence makes it possible to study the effect of a graph
operation on the complex. Consider, for instance, the deletion of an
edge $uv$ from a finite simple graph $G$. The simplices that disappear
from its \defi{clique complex} $\cl(G)$ are precisely those obtained by adjoining
$u$ and $v$ to a clique, possibly empty, of their common neighbourhood.
It is therefore natural to look to this neighbourhood for a condition
under which deleting $uv$ preserves the homotopy type of $\cl(G)$.

Ivashchenko's contractible transformations use such a condition
for edges and an analogous one for vertices
\cite{Ivashchenko1993,Ivashchenko1994a}. Starting from the
one-vertex graph $K_1$, they generate the family $\I$ of \defi{I-contractible
graphs} by allowing vertex and edge additions as well as deletions.
A vertex operation is allowed when the neighbours of the vertex induce
a graph in $\I$. For an edge $uv$, it is the graph $N_G(u,v)$ induced
by the common neighbours that must belong to $\I$. Each transformation
preserves the homotopy type of the clique complex \cite{ChenYauYeh2001}.
Conversely, by a theorem of Dochtermann and Matsushita, every finite
graph with a contractible clique complex belongs to $\I$
\cite[Corollary~1.5, p.~2]{DochtermannMatsushita2025}.

An I-contractible graph can therefore be transformed into $K_1$, but
the sequence of transformations may include additions. To express what deletions alone can
accomplish, we impose the same restriction on the neighbourhood used
at each step. The resulting family $\IS$ of \defi{strong I-contractible
graphs} contains $K_1$. A graph of order at least two belongs to
$\IS$ precisely when it has a vertex or edge for which both the neighbourhood
and the graph left by deletion belong to $\IS$. Such a reduction also yields a simplicial collapse of
$\cl(G)$ to a point \cite[Theorem~1]{EFH2022}.

The distinction between the two families already appears in a flag
triangulation of the dunce hat. Alongside their enumeration of contractible
and dismantlable graphs through order $11$, Dochtermann, Espinoza,
Fr\'ias-Armenta and Hern\'andez \cite{DEFH2023} exhibited such a
triangulation with graph $\Dunce$ on $15$ vertices. It lies in
$\I\setminus\IS$, yet the addition of just one edge gives a graph
$\Dunce^+\in\IS$. The authors conjectured that both graphs had the
least possible number of vertices for their respective properties.
For $\Dunce$ they proposed a stronger minimality assertion using only
the reduced homology of the clique complex over $\mathbb F_2$.
Throughout, the homology of a graph means the reduced homology of its
clique complex, and unqualified acyclicity refers to these coefficients.
We restate the two conjectures in this notation.

\begin{quote}
\textbf{Conjecture 14 of \cite[p.~8]{DEFH2023}.}
Among graphs outside $\IS$ with trivial reduced $\mathbb F_2$-homology,
$\Dunce$ has minimum order. In particular, it has minimum order among
graphs in $\I\setminus\IS$.

\smallskip
\textbf{Conjecture 15 of \cite[p.~9]{DEFH2023}.}
Among graphs $G\in\IS$ having an edge $uv$ such that
$N_G(u,v)\in\IS$ and $G-uv\notin\IS$, the graph $\Dunce^+$ has
minimum order.
\end{quote}

The two minima in Conjecture~14 differ. For acyclicity over
$\mathbb F_2$, the minimum order is $13$. Up to isomorphism, there are exactly three graphs of this
order outside $\IS$ that are acyclic over $\mathbb F_2$. Each has
first integral homology $\ZZ/3$, so changing coefficients to
$\mathbb F_2$ conceals the obstruction to contractibility. The three
examples arise from a single triangulated disc and two successive
edge additions. In conjunction with Adamaszek's classification
\cite{Adamaszek2014}, the disc construction also shows that $13$ is
the least number of vertices of any flag complex with odd torsion in
its first integral homology (Proposition~\ref{rem:adamaszek}).

For contractibility itself the first example occurs one vertex later.
We construct a graph on $14$ vertices whose clique complex is
contractible but has no free face. In the graph, no vertex or edge has
a neighbourhood in $\IS$. Adding one edge nevertheless gives a graph
in $\IS$, and the common neighbours of its endpoints induce a path.
The original graph thus refutes the second assertion of Conjecture~14,
while the one-edge extension refutes Conjecture~15. The example also
attains the least number of vertices of a contractible, non-collapsible
flag complex. Indeed, every contractible flag complex on at most $13$
vertices has its graph in $\IS$ and hence collapses to a point
(Corollary~\ref{cor:minimality}). Earlier examples of contractible,
non-collapsible flag complexes arise from Bing's house
\cite{ChenYauYeh2001} and triangulations of the dunce hat
\cite{FieuxJouve2020,DEFH2023}.

To obtain the lower bounds for the two conjectures, we enumerate all
connected graphs through order $13$. The possible first deletions
determine membership in $\IS$ and in the dismantling levels of Fieux
and Jouve \cite{FieuxJouve2020}. For the acyclic graphs left outside
$\IS$, integral homology or an explicit sequence of transformations
settles I-contractibility. This extends the counts of \cite{DEFH2023}
and gives the minimum order $14$ for both the second assertion of
Conjecture~14 and Conjecture~15. At order $14$ we extend the enumeration to all connected
$K_4$-free graphs. Here every vertex neighbourhood is triangle-free,
so admissibility has a simple local criterion.
Exactly eight such graphs, up to isomorphism, belong to
$\I\setminus\IS$, and the least number of edges among them is $44$.
Section~3 of the accompanying supplement, \emph{Finite closures and graph certificates}, constructs a further $21$ members of $\I\setminus\IS$
of order $14$ that contain a $K_4$. The total number of isomorphism
classes with these latter properties is unknown.

\section{Graph reductions and homology}\label{sec:prelim}

The reductions considered here depend on the neighbourhood of a
vertex or of a pair of vertices. For $v\in V(G)$, let $N_G(v)$ be the subgraph induced by the neighbours
of $v$. For distinct $u,v\in V(G)$, their
\defi{common neighbourhood} $N_G(u,v)$ is the subgraph induced by
their common neighbours. In the clique complex $\cl(G)$,
the link of $v$ is $\cl(N_G(v))$ and the link of an edge $uv$ is
$\cl(N_G(u,v))$. We write $G-v$ and $G-uv$ for the graphs obtained by deleting
the vertex $v$ and the edge $uv$, respectively, and $G+uv$ for the
graph obtained by adding a missing edge $uv$.

\subsection{I-contractibility and its strong variants}
\label{sec:contractible}
Ivashchenko's transformations allow a vertex or an edge to be added
or deleted whenever the appropriate neighbourhood is itself
I-contractible.

\begin{definition}[{\cite[Definition~1]{DEFH2023}}]\label{def:I}
The class $\I$ of \defi{I-contractible graphs} is the least class
containing $K_1$ and closed under the following operations on a graph $G$.
\begin{enumerate}[label=(\roman*),leftmargin=2.5em]
\item delete a vertex $v$ when $N_G(v)\in\I$;
\item add a new vertex whose neighbourhood is an induced subgraph
of $G$ belonging to $\I$;
\item delete an edge $uv$ when $N_G(u,v)\in\I$;
\item add a missing edge $uv$ when $N_G(u,v)\in\I$.
\end{enumerate}
We call these four operations \defi{I-transformations}.
The empty graph does not belong to $\I$.
\end{definition}

Ivashchenko introduced these transformations \cite{Ivashchenko1993}
and showed that they preserve homology \cite{Ivashchenko1994a}.
In fact, each transformation preserves the simple homotopy type of
the clique complex, as Chen, Yau and Yeh showed \cite{ChenYauYeh2001}
(see also \cite[Theorem~2.8, p.~5]{DochtermannMatsushita2025}).
Starting from $K_1$ therefore produces only graphs with contractible
clique complexes. The converse theorem of Dochtermann and Matsushita
identifies $\I$ with all such graphs \cite[Corollary~1.5, p.~2]{DochtermannMatsushita2025}.

For a reduction to $K_1$ by deletions alone, the same restriction can
be imposed recursively on the neighbourhood at each step. This gives
the following two subfamilies of $\I$.

\begin{definition}[{\cite[Definition~3 and Lemma~1]{EFH2022};
\cite[Definitions~4--5]{DEFH2023}}]\label{def:IS}
Let $\IS$ denote the class of \defi{strong I-contractible graphs}.
It contains $K_1$. For every graph $G$ of order at least two,
membership in $\IS$ is determined recursively by
\[
 G\in\IS\quad\Longleftrightarrow\quad
 \begin{cases}
 G-v,\ N_G(v)\in\IS&\text{for some }v\in V(G),\quad\text{or}\\
 G-uv,\ N_G(u,v)\in\IS&\text{for some }uv\in E(G).
 \end{cases}
\]
The class $\ISV$ of \defi{strong vertex I-contractible graphs} also
contains $K_1$. For a graph $G$ of order at least two, $G\in\ISV$ if
and only if there is a vertex $v$ for which both $G-v$ and $N_G(v)$
belong to $\ISV$. The empty graph belongs to neither class.
\end{definition}

A deletion decreases $|V|+|E|$, and the neighbourhood tested at that
step also has smaller $|V|+|E|$ than $G$. Induction on this quantity
therefore gives
\begin{equation}\label{eq:strong-inclusions}
 \ISV\subseteq\IS\subseteq\I.
\end{equation}
For the first inclusion, every vertex deletion allowed by the
recursion for $\ISV$ is also allowed in $\IS$, by the induction
hypothesis applied to the neighbourhood and the remaining graph.
For the second, take a deletion supplied by Definition~\ref{def:IS}.
Both of these smaller graphs now belong to $\I$ by induction, so
operation (ii) or (iv) of Definition~\ref{def:I} recovers $G$.
Thus reversing a reduction to $K_1$ gives a construction from $K_1$
by the corresponding additions.

We call a vertex $v$ \defi{$\IS$-admissible} if $N_G(v)\in\IS$,
and a pair of distinct vertices $u,v$ \defi{$\IS$-admissible} if
$N_G(u,v)\in\IS$. We use \defi{admissible} in this sense throughout.
An admissible vertex can be deleted, and a new vertex can be added
when its neighbourhood in the resulting graph belongs to $\IS$.
An edge is admissible if its endpoints form an admissible pair.
We may delete such an edge or add a missing edge between the vertices
of an admissible pair. A \defi{strong reduction} is a sequence of
admissible deletions ending at $K_1$. Since $\emptyset\notin\IS$,
an admissible deletion cannot remove an isolated vertex or an edge
whose endpoints have no common neighbour.

On triangle-free graphs these conditions have a particularly simple
form. Every nonempty tree belongs to $\ISV$, since a leaf has
neighbourhood $K_1$ and successive leaf deletions leave a single
vertex. Conversely, the clique complex of a triangle-free graph is
the graph itself. If the graph belongs to $\I$, it is connected
and has trivial first homology, so it is a tree.

\subsection{Dismantling levels}
The condition on the neighbourhood of a vertex can also be imposed
one level at a time. In the dismantling hierarchy, a reduction at
level $k$ uses vertex deletions whose neighbourhoods belong to the
preceding level.

\begin{definition}[{\cite[Definition~1]{FieuxJouve2020}}]\label{def:kdis}
A \defi{cone} is a nonempty graph with a vertex adjacent to every other
vertex. Cones form the class $(-1)$-Dis. For $k\geq0$, a vertex $v$
is \defi{$k$-dismantlable} when $N_G(v)$ belongs to $(k-1)$-Dis.
A graph is \defi{$k$-dismantlable}, or belongs to $k$-Dis, if it can
be reduced to $K_1$ by successive deletions of $k$-dismantlable
vertices. The graph $K_1$ belongs to every level. The empty graph belongs to none.
\end{definition}

At level zero the neighbourhood of a vertex $v$ must be a cone.
Its apex is a neighbour $u$ of $v$ adjacent to every other neighbour
of $v$. Equivalently, $N_G[v]\subseteq N_G[u]$ for some $u\ne v$,
where $N_G[v]=\{v\}\cup V(N_G(v))$ denotes the closed neighbourhood
of $v$.

Every reduction at one level is also available at the next \cite[Proposition~2]{FieuxJouve2020}. Together with
\cite[Lemma~19 and Theorem~23]{DEFH2023}, this gives
\begin{equation}\label{eq:dismantling-inclusions}
 (-1)\text{-Dis}\subseteq0\text{-Dis}\subseteq1\text{-Dis}
 \subseteq\cdots\subseteq\ISV\subseteq\IS\subseteq\I.
\end{equation}
In fact, every graph in $\ISV$ occurs at some finite level
\cite[Theorem~23]{DEFH2023}.

\subsection{Homology and acyclicity}
Membership in $\I$ forces the reduced homology of the clique complex
to vanish. To compare this condition for different coefficients, write $\widetilde H_i(G;R)=\widetilde H_i(\cl(G);R)$ and,
for $R=\ZZ$ or $\mathbb F_2$, let $\mathcal A_R$ denote the
nonempty graphs whose reduced homology with coefficients in $R$
vanishes in every degree $i\geq0$. A graph in either family has
Euler characteristic one. This characteristic can be read directly
from the graph, since
\[
 \chi(G)=\sum_{i\geq0}(-1)^if_i(G),
\]
where $f_i(G)$ counts the cliques with $i+1$ vertices.

The passage from integral homology to homology over $\mathbb F_2$
is governed by the universal coefficient theorem
\cite[Theorem~3A.3]{Hatcher2002}. For a nonempty finite complex $K$
and $i\geq0$, it gives
\begin{equation}\label{eq:uct}
 0\longrightarrow\widetilde H_i(K;\ZZ)\otimes\mathbb F_2
 \longrightarrow\widetilde H_i(K;\mathbb F_2)
 \longrightarrow\operatorname{Tor}^{\ZZ}_1
       (\widetilde H_{i-1}(K;\ZZ),\mathbb F_2)
 \longrightarrow0.
\end{equation}
Here $\widetilde H_{-1}(K;\ZZ)=0$ because $K$ is nonempty.
For an integrally acyclic complex both outer terms vanish in every
degree, so its reduced homology over $\mathbb F_2$ also vanishes.
Since every contractible complex is integrally acyclic, we obtain
\begin{equation}\label{eq:homology-inclusions}
 \IS\subseteq\I\subseteq\mathcal A_{\ZZ}
 \subseteq\mathcal A_{\mathbb F_2}.
\end{equation}
Passing to $\mathbb F_2$ loses information about odd torsion. In
particular, a summand $\ZZ/3$ in integral homology contributes
nothing to either outer term of~\eqref{eq:uct}, since both
$(\ZZ/3)\otimes\mathbb F_2$ and
$\operatorname{Tor}^{\ZZ}_1(\ZZ/3,\mathbb F_2)$ vanish.

\subsection{A disc with boundary identifications}\label{sec:discs}
A disc attached to a circle by a map of degree $3$ gives the Moore
space $M(\ZZ/3,1)$. More generally, if the boundary is attached to
a connected graph with a single cycle, contracting a maximal tree
reduces the attaching map to a map between circles. Its degree then
determines the homotopy type.

\begin{proposition}\label{prop:disc}
Let $q:D\to X$ be a simplicial map from a triangulated closed disc.
Form a quotient of $D$ by identifying boundary points with the same
image under $q$ and making no other identifications. Suppose that
$q$ induces a homeomorphism from this quotient onto $X$. If
$S=q(\partial D)$ is a connected graph with exactly one cycle, then
\[
 X\simeq S^1\cup_d D^2,
\]
where $d$ is the degree of the boundary map after contracting a maximal
tree in $S$. For $|d|=1$ the complex is contractible. For $|d|=3$
it is homotopy equivalent to $M(\ZZ/3,1)$ and has
$\widetilde H_1(X;\ZZ)\cong\ZZ/3$ and zero reduced homology in all
other degrees.
\end{proposition}
\begin{proof}
Write $f=q|_{\partial D}:\partial D\to S$. The hypothesis on $q$
means that $X$ is obtained by attaching $D$ to $S$ along $f$.
Choose a maximal tree $T$ in $S$. Contracting $T$ preserves the
homotopy type of $X$ because $T$ is a contractible subcomplex
\cite[Propositions~0.16--0.17]{Hatcher2002}. With $p:S\to S/T$
the quotient map, the resulting space is
$X/T=(S/T)\cup_{pf}D^2$. Exactly one edge of $S$ lies outside $T$,
so $S/T\cong S^1$.

The degree identifies $\pi_1(S^1)$ with $\ZZ$
\cite[Theorem~1.7]{Hatcher2002}. Its conjugacy classes are singletons,
so the degree also determines the free homotopy class of $pf$.
Since homotopic attaching maps give homotopy equivalent spaces
\cite[Proposition~0.18]{Hatcher2002}, this proves the asserted
homotopy equivalence. Reflecting $D$ reverses the sign of the degree
and gives a homeomorphic adjunction space. When $|d|=1$, we may thus
attach the disc by the identity, obtaining a contractible space.
When $|d|=3$, we obtain the stated Moore space. Its cellular chain
complex has one generator in each of dimensions zero, one and two,
and the boundary from degree two to degree one is multiplication by
$d$ on $\ZZ$. This gives the asserted reduced homology.
\end{proof}

\section{Counterexamples on thirteen and fourteen vertices}\label{sec:examples}

The two assertions of Conjecture~14 lead to different constructions.
On thirteen vertices we realise odd torsion in a clique complex whose
reduced homology over $\mathbb F_2$ vanishes. This gives three graphs
outside $\IS$ that refute the first assertion. To refute the
I-contractibility assertion, we construct three graphs on fourteen
vertices that belong to $\I$ but admit no first admissible deletion.
These graphs also refute the first assertion. In each case, adding
one edge gives a graph in $\IS$ that refutes Conjecture~15.

The adjacency lists in Table~\ref{tab:witness-adjacency} specify the
initial graph $\Gt$ on $\{0,\ldots,12\}$ and the three graphs
$\Ga,\Gb,\Gc$ on $\{0,\ldots,13\}$.

\begin{table}[htbp]
\centering\small
\caption{Adjacency lists for $\Gt,\Ga,\Gb,\Gc$. Row $i$ lists the
neighbours of $i$ with labels greater than $i$. An empty set means that $i$ has
no such neighbour, and a dash that $i$ is absent.}\label{tab:witness-adjacency}

\begin{tabular}{@{}rllll@{}}
\toprule
$i$ & $\Gt$ & $\Ga$ & $\Gb$ & $\Gc$ \\
\midrule
0 & $4,5,9,10,11,12$ & $4,5,8,10,12,13$ & $5,7,9,12,13$ & $4,5,7,9,12,13$ \\
1 & $4,6,8,9,11,12$ & $5,7,8,10,12,13$ & $5,8,9,11,13$ & $4,5,8,10,12,13$ \\
2 & $5,7,8,9,10,12$ & $6,7,9,12,13$ & $6,7,8,11,12,13$ & $6,7,8,10,11,13$ \\
3 & $6,7,8,10,11,12$ & $6,9,10,12,13$ & $6,7,9,10,12,13$ & $6,7,9,10,11,13$ \\
4 & $8,10,12$ & $8,9,11,13$ & $6,8,9,10,11,13$ & $8,9,11,13$ \\
5 & $8,11,12$ & $9,10,11,12$ & $7,8,10,11,12$ & $9,10,11,13$ \\
6 & $9,10,12$ & $7,10,11,12$ & $10,11,12$ & $8,9,12,13$ \\
7 & $9,11,12$ & $8,11,13$ & $10,11,13$ & $10,12,13$ \\
8 & $10,11$ & $11,12$ & $10,12,13$ & $11,12$ \\
9 & $10,11$ & $11,12,13$ & $11,12,13$ & $11,12$ \\
10 & $\varnothing$ & $11,13$ & $\varnothing$ & $11,12$ \\
11 & $\varnothing$ & $\varnothing$ & $\varnothing$ & $\varnothing$ \\
12 & $\varnothing$ & $\varnothing$ & $\varnothing$ & $\varnothing$ \\
13 & --- & $\varnothing$ & $\varnothing$ & $\varnothing$ \\
\bottomrule
\end{tabular}

\end{table}

\subsection{Odd torsion on thirteen vertices}\label{sec:thirteen}

In $\cl(\Gt)$, four edges form the circle to which a triangulated
disc is attached. The attaching map has degree three, and two
successive edge additions preserve the resulting homotopy type.

\begin{theorem}\label{thm:A}
Let $\Gt'=\Gt+\{4,7\}$ and $\Gt''=\Gt'+\{5,6\}$.
The graphs $\Gt,\Gt',\Gt''$ have $40$, $41$ and $42$ edges, respectively,
and their clique complexes are homotopy equivalent to $M(\ZZ/3,1)$.
In particular, they are $\Ztwo$-acyclic and none belongs to $\I$.
They refute the first assertion of Conjecture~14.
\end{theorem}

\begin{figure}[tbp]
\centering
\begin{tikzpicture}[x=1cm,y=1cm]
\definecolor{witnessred}{RGB}{181,55,36}
\definecolor{witnessblue}{RGB}{32,89,140}
\coordinate (v0) at (0.942769700,-1.632925019);
\coordinate (v1) at (1.632925019,0.942769700);
\coordinate (v2) at (-1.632925019,-0.942769700);
\coordinate (v3) at (-0.942769700,1.632925019);
\coordinate (v4) at (1.981676324,-0.530988571);
\coordinate (v5) at (-0.530988571,-1.981676324);
\coordinate (v6) at (0.530988571,1.981676324);
\coordinate (v7) at (-1.981676324,0.530988571);
\coordinate (v8) at (-1.700000000,2.944486373);
\coordinate (v9) at (3.400000000,0.000000000);
\coordinate (v10) at (-1.700000000,-2.944486373);
\coordinate (v11) at (1.700000000,-2.944486373);
\coordinate (v12) at (1.700000000,2.944486373);
\coordinate (v13) at (-3.400000000,-0.000000000);
\coordinate (v14) at (-2.944486373,-1.700000000);
\coordinate (v15) at (2.944486373,-1.700000000);
\coordinate (v16) at (0.000000000,3.400000000);
\coordinate (v17) at (2.944486373,1.700000000);
\coordinate (v18) at (0.000000000,-3.400000000);
\coordinate (v19) at (-2.944486373,1.700000000);
\coordinate (v20) at (0.000000000,0.000000000);
\fill[blue!3] (v0) -- (v4) -- (v15) -- cycle;
\fill[blue!3] (v0) -- (v4) -- (v20) -- cycle;
\fill[blue!3] (v0) -- (v5) -- (v18) -- cycle;
\fill[blue!3] (v0) -- (v5) -- (v20) -- cycle;
\fill[blue!3] (v0) -- (v11) -- (v15) -- cycle;
\fill[blue!3] (v0) -- (v11) -- (v18) -- cycle;
\fill[blue!3] (v1) -- (v4) -- (v9) -- cycle;
\fill[blue!3] (v1) -- (v4) -- (v20) -- cycle;
\fill[blue!3] (v1) -- (v6) -- (v12) -- cycle;
\fill[blue!3] (v1) -- (v6) -- (v20) -- cycle;
\fill[blue!3] (v1) -- (v9) -- (v17) -- cycle;
\fill[blue!3] (v1) -- (v12) -- (v17) -- cycle;
\fill[blue!3] (v2) -- (v5) -- (v10) -- cycle;
\fill[blue!3] (v2) -- (v5) -- (v20) -- cycle;
\fill[blue!3] (v2) -- (v7) -- (v13) -- cycle;
\fill[blue!3] (v2) -- (v7) -- (v20) -- cycle;
\fill[blue!3] (v2) -- (v10) -- (v14) -- cycle;
\fill[blue!3] (v2) -- (v13) -- (v14) -- cycle;
\fill[blue!3] (v3) -- (v6) -- (v16) -- cycle;
\fill[blue!3] (v3) -- (v6) -- (v20) -- cycle;
\fill[blue!3] (v3) -- (v7) -- (v19) -- cycle;
\fill[blue!3] (v3) -- (v7) -- (v20) -- cycle;
\fill[blue!3] (v3) -- (v8) -- (v16) -- cycle;
\fill[blue!3] (v3) -- (v8) -- (v19) -- cycle;
\fill[blue!3] (v4) -- (v9) -- (v15) -- cycle;
\fill[blue!3] (v5) -- (v10) -- (v18) -- cycle;
\fill[blue!3] (v6) -- (v12) -- (v16) -- cycle;
\fill[blue!3] (v7) -- (v13) -- (v19) -- cycle;
\draw[draw=black!52,line width=.38pt] (v0) -- (v4);
\draw[draw=black!52,line width=.38pt] (v0) -- (v5);
\draw[draw=black!52,line width=.38pt] (v0) -- (v11);
\draw[draw=black!52,line width=.38pt] (v0) -- (v15);
\draw[draw=black!52,line width=.38pt] (v0) -- (v18);
\draw[draw=black!52,line width=.38pt] (v0) -- (v20);
\draw[draw=black!52,line width=.38pt] (v1) -- (v4);
\draw[draw=black!52,line width=.38pt] (v1) -- (v6);
\draw[draw=black!52,line width=.38pt] (v1) -- (v9);
\draw[draw=black!52,line width=.38pt] (v1) -- (v12);
\draw[draw=black!52,line width=.38pt] (v1) -- (v17);
\draw[draw=black!52,line width=.38pt] (v1) -- (v20);
\draw[draw=black!52,line width=.38pt] (v2) -- (v5);
\draw[draw=black!52,line width=.38pt] (v2) -- (v7);
\draw[draw=black!52,line width=.38pt] (v2) -- (v10);
\draw[draw=black!52,line width=.38pt] (v2) -- (v13);
\draw[draw=black!52,line width=.38pt] (v2) -- (v14);
\draw[draw=black!52,line width=.38pt] (v2) -- (v20);
\draw[draw=black!52,line width=.38pt] (v3) -- (v6);
\draw[draw=black!52,line width=.38pt] (v3) -- (v7);
\draw[draw=black!52,line width=.38pt] (v3) -- (v8);
\draw[draw=black!52,line width=.38pt] (v3) -- (v16);
\draw[draw=black!52,line width=.38pt] (v3) -- (v19);
\draw[draw=black!52,line width=.38pt] (v3) -- (v20);
\draw[draw=black!52,line width=.38pt] (v4) -- (v9);
\draw[draw=black!52,line width=.38pt] (v4) -- (v15);
\draw[draw=black!52,line width=.38pt] (v4) -- (v20);
\draw[draw=black!52,line width=.38pt] (v5) -- (v10);
\draw[draw=black!52,line width=.38pt] (v5) -- (v18);
\draw[draw=black!52,line width=.38pt] (v5) -- (v20);
\draw[draw=black!52,line width=.38pt] (v6) -- (v12);
\draw[draw=black!52,line width=.38pt] (v6) -- (v16);
\draw[draw=black!52,line width=.38pt] (v6) -- (v20);
\draw[draw=black!52,line width=.38pt] (v7) -- (v13);
\draw[draw=black!52,line width=.38pt] (v7) -- (v19);
\draw[draw=black!52,line width=.38pt] (v7) -- (v20);
\draw[draw=witnessred,line width=.8pt] (v8) -- (v16);
\draw[draw=witnessred,line width=.8pt] (v8) -- (v19);
\draw[draw=witnessred,line width=.8pt] (v9) -- (v15);
\draw[draw=witnessred,line width=.8pt] (v9) -- (v17);
\draw[draw=witnessred,line width=.8pt] (v10) -- (v14);
\draw[draw=witnessred,line width=.8pt] (v10) -- (v18);
\draw[draw=witnessred,line width=.8pt] (v11) -- (v15);
\draw[draw=witnessred,line width=.8pt] (v11) -- (v18);
\draw[draw=witnessred,line width=.8pt] (v12) -- (v16);
\draw[draw=witnessred,line width=.8pt] (v12) -- (v17);
\draw[draw=witnessred,line width=.8pt] (v13) -- (v14);
\draw[draw=witnessred,line width=.8pt] (v13) -- (v19);
\node[jfef vertex] at (v0) {$0$};
\node[jfef vertex] at (v1) {$1$};
\node[jfef vertex] at (v2) {$2$};
\node[jfef vertex] at (v3) {$3$};
\node[jfef vertex] at (v4) {$4$};
\node[jfef vertex] at (v5) {$5$};
\node[jfef vertex] at (v6) {$6$};
\node[jfef vertex] at (v7) {$7$};
\node[jfef vertex,draw=witnessred,line width=.7pt] at (v8) {$8$};
\node[jfef vertex,draw=witnessred,line width=.7pt] at (v9) {$8$};
\node[jfef vertex,draw=witnessred,line width=.7pt] at (v10) {$8$};
\node[jfef vertex,draw=witnessred,line width=.7pt] at (v11) {$9$};
\node[jfef vertex,draw=witnessred,line width=.7pt] at (v12) {$9$};
\node[jfef vertex,draw=witnessred,line width=.7pt] at (v13) {$9$};
\node[jfef vertex,draw=witnessred,line width=.7pt] at (v14) {$10$};
\node[jfef vertex,draw=witnessred,line width=.7pt] at (v15) {$10$};
\node[jfef vertex,draw=witnessred,line width=.7pt] at (v16) {$10$};
\node[jfef vertex,draw=witnessred,line width=.7pt] at (v17) {$11$};
\node[jfef vertex,draw=witnessred,line width=.7pt] at (v18) {$11$};
\node[jfef vertex,draw=witnessred,line width=.7pt] at (v19) {$11$};
\node[jfef vertex] at (v20) {$12$};
\end{tikzpicture}\hspace{6mm}\begin{tikzpicture}[x=1cm,y=1cm]
\definecolor{witnessred}{RGB}{181,55,36}
\definecolor{witnessblue}{RGB}{32,89,140}
\coordinate (v0) at (1.146661411,2.139452559);
\coordinate (v1) at (-2.265488516,1.039851549);
\coordinate (v2) at (2.265489343,-1.039852163);
\coordinate (v3) at (-1.146663831,-2.139454358);
\coordinate (v4) at (-1.146661411,2.139452559);
\coordinate (v5) at (2.265488516,1.039851549);
\coordinate (v6) at (-2.265489343,-1.039852163);
\coordinate (v7) at (1.146663831,-2.139454358);
\coordinate (v8) at (3.440000000,-3.200000000);
\coordinate (v9) at (-3.440000000,3.200000000);
\coordinate (v10) at (3.440000000,3.200000000);
\coordinate (v11) at (-3.440000000,-3.200000000);
\coordinate (v12) at (0.000000000,0.000000000);
\draw[draw=black!52,line width=.38pt] (v0) -- (v4);
\draw[draw=black!52,line width=.38pt] (v0) -- (v5);
\draw[draw=black!52,line width=.38pt] (v0) -- (v9);
\draw[draw=black!52,line width=.38pt] (v0) -- (v10);
\draw[draw=black!52,line width=.38pt] (v0) -- (v11);
\draw[draw=black!52,line width=.38pt] (v0) -- (v12);
\draw[draw=black!52,line width=.38pt] (v1) -- (v4);
\draw[draw=black!52,line width=.38pt] (v1) -- (v6);
\draw[draw=black!52,line width=.38pt] (v1) -- (v8);
\draw[draw=black!52,line width=.38pt] (v1) -- (v9);
\draw[draw=black!52,line width=.38pt] (v1) -- (v11);
\draw[draw=black!52,line width=.38pt] (v1) -- (v12);
\draw[draw=black!52,line width=.38pt] (v2) -- (v5);
\draw[draw=black!52,line width=.38pt] (v2) -- (v7);
\draw[draw=black!52,line width=.38pt] (v2) -- (v8);
\draw[draw=black!52,line width=.38pt] (v2) -- (v9);
\draw[draw=black!52,line width=.38pt] (v2) -- (v10);
\draw[draw=black!52,line width=.38pt] (v2) -- (v12);
\draw[draw=black!52,line width=.38pt] (v3) -- (v6);
\draw[draw=black!52,line width=.38pt] (v3) -- (v7);
\draw[draw=black!52,line width=.38pt] (v3) -- (v8);
\draw[draw=black!52,line width=.38pt] (v3) -- (v10);
\draw[draw=black!52,line width=.38pt] (v3) -- (v11);
\draw[draw=black!52,line width=.38pt] (v3) -- (v12);
\draw[draw=black!52,line width=.38pt] (v4) -- (v8);
\draw[draw=black!52,line width=.38pt] (v4) -- (v10);
\draw[draw=black!52,line width=.38pt] (v4) -- (v12);
\draw[draw=black!52,line width=.38pt] (v5) -- (v8);
\draw[draw=black!52,line width=.38pt] (v5) -- (v11);
\draw[draw=black!52,line width=.38pt] (v5) -- (v12);
\draw[draw=black!52,line width=.38pt] (v6) -- (v9);
\draw[draw=black!52,line width=.38pt] (v6) -- (v10);
\draw[draw=black!52,line width=.38pt] (v6) -- (v12);
\draw[draw=black!52,line width=.38pt] (v7) -- (v9);
\draw[draw=black!52,line width=.38pt] (v7) -- (v11);
\draw[draw=black!52,line width=.38pt] (v7) -- (v12);
\draw[draw=witnessred,line width=.8pt] (v8) -- (v10);
\draw[draw=witnessred,line width=.8pt] (v8) -- (v11);
\draw[draw=witnessred,line width=.8pt] (v9) -- (v10);
\draw[draw=witnessred,line width=.8pt] (v9) -- (v11);
\node[jfef vertex] at (v0) {$0$};
\node[jfef vertex] at (v1) {$1$};
\node[jfef vertex] at (v2) {$2$};
\node[jfef vertex] at (v3) {$3$};
\node[jfef vertex] at (v4) {$4$};
\node[jfef vertex] at (v5) {$5$};
\node[jfef vertex] at (v6) {$6$};
\node[jfef vertex] at (v7) {$7$};
\node[jfef vertex,draw=witnessred,line width=.7pt] at (v8) {$8$};
\node[jfef vertex,draw=witnessred,line width=.7pt] at (v9) {$9$};
\node[jfef vertex,draw=witnessred,line width=.7pt] at (v10) {$10$};
\node[jfef vertex,draw=witnessred,line width=.7pt] at (v11) {$11$};
\node[jfef vertex] at (v12) {$12$};
\end{tikzpicture}
\caption{The triangulated disc with $28$ triangles (left) and the graph
$\Gt$ (right). Read clockwise from the uppermost vertex, labelled
$10$, the boundary follows the highlighted cycle
$10$--$9$--$11$--$8$--$10$ three times. The starting vertex is
$b_0$ in \ref{app:disc-certificates}.}\label{fig:G13}
\end{figure}
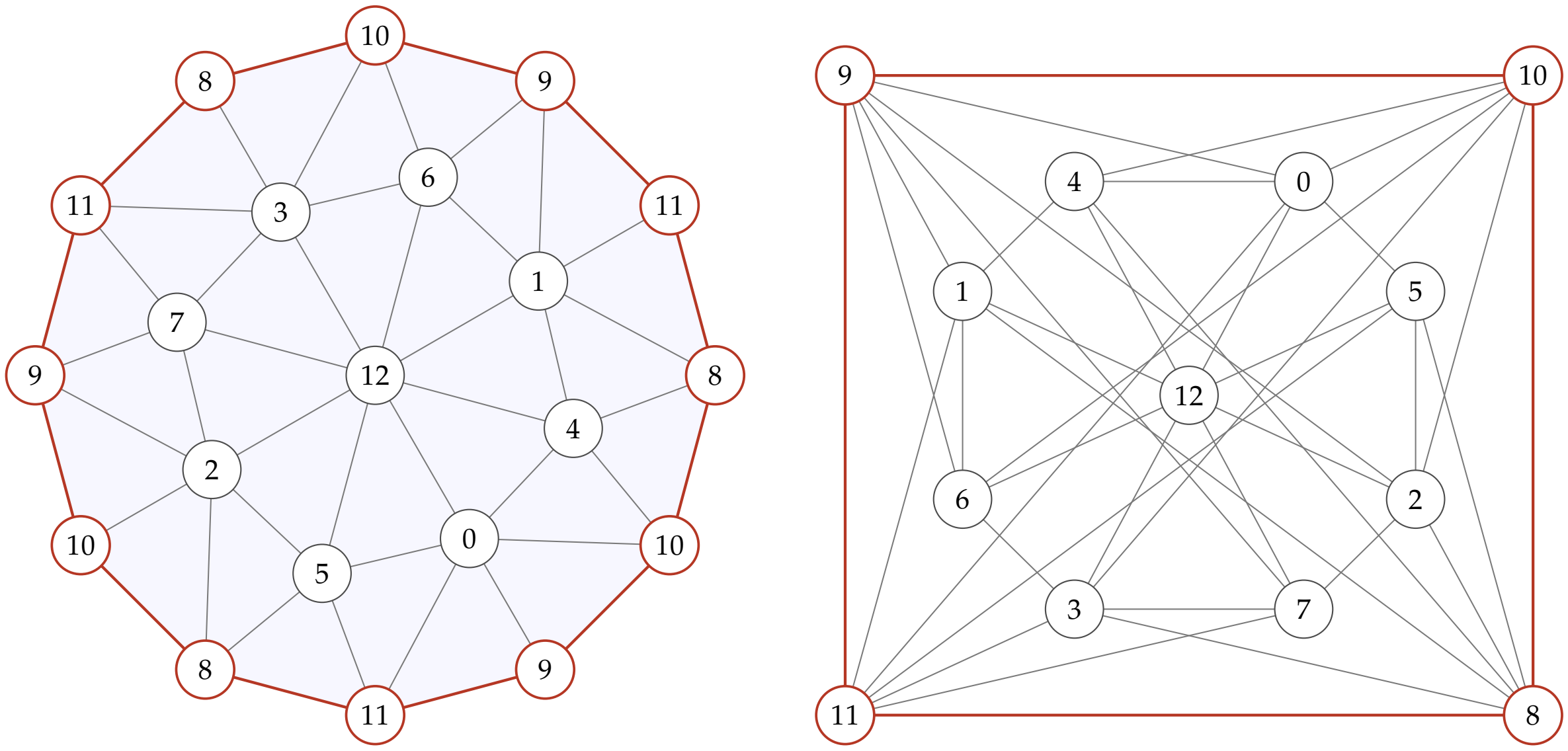

\begin{proof}
The adjacency lists show that $\Gt$ is $K_4$-free and has $28$ triangles.
The four edges of the cycle $10\text{--}9\text{--}11\text{--}8\text{--}10$
each lie in three triangles, and every other edge lies in two.
To see how these triangles meet along the cycle, consider the disc
on the left of Figure~\ref{fig:G13}. Its triangles are the ones
listed in Table~\ref{tab:G13-disc-triangles}, attached successively
along boundary arcs. Reading clockwise from the top gives the boundary
labels
\[
 10,9,11,8\mid10,9,11,8\mid10,9,11,8.
\]
Each interior vertex has its own label, distinct from all boundary
labels. Identifying boundary vertices with the same label and
the corresponding boundary edges produces $\cl(\Gt)$
(\ref{app:disc-certificates}). The three occurrences of
the same cyclic word give the attaching map degree three. Thus
Proposition~\ref{prop:disc} gives
\[
 \cl(\Gt)\simeq M(\ZZ/3,1),\qquad
 \widetilde H_1(\Gt;\ZZ)\cong\ZZ/3,
 \qquad \widetilde H_i(\Gt;\ZZ)=0\quad(i\ne1).
\]
The only nonzero reduced integral homology group has order three,
so~\eqref{eq:uct} gives $\widetilde H_i(\Gt;\Ztwo)=0$ in every degree.

For the addition of $\{4,7\}$ to $\Gt$ and then of $\{5,6\}$ to $\Gt'$,
the unique common neighbour of the endpoints is vertex $12$ in the
graph at that step. Both additions are therefore admissible and
preserve the homotopy type of the clique complex. All three graphs
are consequently $\mathbb F_2$-acyclic, while their nonzero integral
homology excludes them from $\I$ and hence from $\IS$. Each has fewer
than the $15$ vertices proposed in the first assertion of Conjecture~14.
\end{proof}

The order of $\Gt$ is already forced by the odd torsion in its first
homology. The lower bound follows from Adamaszek's classification
together with the effect of adjoining one vertex.

\begin{proposition}\label{rem:adamaszek}
A flag complex with odd torsion in its first integral homology has at
least $13$ vertices, and this bound is attained.
\end{proposition}

\begin{proof}
As in Adamaszek's classification, call $G$ a torsion graph if
$H_1(\cl(G);\ZZ)$ has nonzero torsion, and call it irreducible if
deleting any vertex gives a graph with torsion-free first homology.
By Adamaszek's theorem \cite[Theorem~1.1]{Adamaszek2014}, there is no
torsion graph on at most ten vertices. On eleven vertices
there are exactly four, and their clique complexes are homotopy
equivalent to the real projective plane. Their first homology is
therefore $\ZZ/2$. On twelve vertices there are $363$ irreducible
torsion graphs, each with torsion subgroup $\ZZ/2$ in its first
integral homology \cite[Section~4, p.~229]{Adamaszek2014}. The classification through
order eleven also follows from Katzman's results on independence
complexes \cite[Theorem~4.1 and Section~5]{Katzman2006}, applied to
the complement of $G$.

The remaining possibility is a reducible torsion graph $G$ on twelve
vertices. For some vertex $v$, the graph $G-v$ must then be one of the
four eleven-vertex examples. An isolated $v$ leaves first homology
unchanged, giving $H_1(\cl(G);\ZZ)=H_1(\cl(G-v);\ZZ)\cong\ZZ/2$.
If $v$ is not isolated, its link is nonempty. The Mayer--Vietoris
sequence for the decomposition
\[
 \cl(G)=\cl(G-v)\cup\bigl(v*\cl(\Nb_G(v))\bigr)
\]
gives a short exact sequence
$0\to Q\to H_1(\cl(G);\ZZ)\to F\to0$, where $Q$ is a quotient of
$\ZZ/2$ and $F$ is a subgroup of
$\widetilde H_0(\cl(\Nb_G(v));\ZZ)$. The group $F$ is free abelian,
so the sequence splits. Thus $H_1(\cl(G);\ZZ)$ has no odd torsion.
The graph $\Gt$ attains the bound by Theorem~\ref{thm:A}.
\end{proof}

\subsection{I-contractible graphs on fourteen vertices}
\label{sec:fourteen}

The graphs on thirteen vertices lie outside $\I$ because of their
integral homology. A different obstruction is needed for an example
in $\I\setminus\IS$. In $\Ga$ every vertex and edge fails the
neighbourhood condition for a first deletion, yet adding a single
edge makes a strong reduction possible.

\begin{theorem}\label{thm:B}
The graph $\Ga$ is $K_4$-free, has $14$ vertices and $44$ edges, and
belongs to $\I\setminus\IS$. The graph $\Ga+\{0,1\}$ belongs to
$\IS$, and the deletion of $\{0,1\}$ from this graph is admissible.
Consequently, $\Ga$ refutes the second assertion of Conjecture~14,
and $\Ga+\{0,1\}$ refutes Conjecture~15.
\end{theorem}

\begin{proof}
In $\Ga$ the endpoints of the missing edge $\{0,1\}$ have as their
common neighbourhood the induced path
\[
 8\text{--}12\text{--}5\text{--}10\text{--}13.
\]
Since nonempty trees belong to $\IS$ (Section~\ref{sec:contractible}),
the edge can be added admissibly. We first reduce the resulting
graph to $K_1$. Write $e(u,w)$ for deletion of the edge $\{u,w\}$
and $v(u)$ for deletion of vertex $u$, retaining the labels of the
surviving vertices. A reduction to vertex $13$ is given by
\[
\begin{gathered}
 e(0,5),\ e(0,10),\ e(0,12),\ e(8,12),\ e(1,12),\\
 v(5),\ v(10),\ v(3),\ v(6),\ v(2),\ v(7),\ v(1),\ v(0),\\
 v(4),\ v(8),\ v(11),\ v(12),\ v(9).
\end{gathered}
\]
Table~\ref{tab:Ga-reduction} gives the neighbourhood in the current
graph at every step. Each is a nonempty tree, so this sequence proves
$\Ga+\{0,1\}\in\IS$. Adding $\{0,1\}$ leaves its common neighbourhood
unchanged. We may therefore delete that edge by an I-transformation
and conclude that $\Ga\in\I$.

It remains to exclude a strong reduction of $\Ga$ itself.
Each vertex neighbourhood is a connected triangle-free graph with
positive first Betti number (Table~\ref{tab:witness-links}).
The endpoints of each edge of $\Ga$ have two or three common
neighbours. No two of them can be adjacent, since they would form a
$K_4$ with the endpoints. The vertex neighbourhoods thus have nonzero
first homology and the edge neighbourhoods have nonzero reduced
zeroth homology. By~\eqref{eq:homology-inclusions}, none belongs to
$\IS$. This rules out every possible first deletion. At order $14$
we have therefore obtained a counterexample $\Ga$ to the second
assertion of Conjecture~14 and a counterexample $\Ga+\{0,1\}$ to
Conjecture~15.
\end{proof}

\begin{figure}[tbp]
\centering
\begin{tikzpicture}[x=1cm,y=1cm]
\definecolor{witnessred}{RGB}{181,55,36}
\definecolor{witnessblue}{RGB}{32,89,140}
\coordinate (v0) at (4.899000000,-0.124598525);
\coordinate (v1) at (2.722969139,3.944383940);
\coordinate (v2) at (-2.435697872,0.406103749);
\coordinate (v3) at (-4.899000000,4.045228676);
\coordinate (v4) at (1.271029059,-0.654558239);
\coordinate (v5) at (4.387943889,-3.466555902);
\coordinate (v6) at (-4.899000000,-4.554000000);
\coordinate (v7) at (-2.939555778,-3.846959139);
\coordinate (v8) at (-0.868699866,3.253825943);
\coordinate (v9) at (2.353335221,-2.606664647);
\coordinate (v10) at (3.219437385,-4.554000000);
\coordinate (v11) at (3.996206901,3.171731296);
\coordinate (v12) at (1.198003432,4.416000000);
\coordinate (v13) at (-1.343328314,-3.846959139);
\draw[draw=witnessblue,line width=1pt,dashed] (v0) -- (v1);
\draw[draw=black!52,line width=.38pt] (v0) -- (v4);
\draw[draw=black!52,line width=.38pt] (v0) -- (v5);
\draw[draw=black!52,line width=.38pt] (v0) -- (v8);
\draw[draw=black!52,line width=.38pt] (v0) -- (v10);
\draw[draw=black!52,line width=.38pt] (v0) -- (v12);
\draw[draw=black!52,line width=.38pt] (v0) -- (v13);
\draw[draw=black!52,line width=.38pt] (v1) -- (v5);
\draw[draw=black!52,line width=.38pt] (v1) -- (v7);
\draw[draw=black!52,line width=.38pt] (v1) -- (v8);
\draw[draw=black!52,line width=.38pt] (v1) -- (v10);
\draw[draw=black!52,line width=.38pt] (v1) -- (v12);
\draw[draw=black!52,line width=.38pt] (v1) -- (v13);
\draw[draw=black!52,line width=.38pt] (v2) -- (v6);
\draw[draw=black!52,line width=.38pt] (v2) -- (v7);
\draw[draw=black!52,line width=.38pt] (v2) -- (v9);
\draw[draw=black!52,line width=.38pt] (v2) -- (v12);
\draw[draw=black!52,line width=.38pt] (v2) -- (v13);
\draw[draw=black!52,line width=.38pt] (v3) -- (v6);
\draw[draw=black!52,line width=.38pt] (v3) -- (v9);
\draw[draw=black!52,line width=.38pt] (v3) -- (v10);
\draw[draw=black!52,line width=.38pt] (v3) -- (v12);
\draw[draw=black!52,line width=.38pt] (v3) -- (v13);
\draw[draw=black!52,line width=.38pt] (v4) -- (v8);
\draw[draw=black!52,line width=.38pt] (v4) -- (v9);
\draw[draw=black!52,line width=.38pt] (v4) -- (v11);
\draw[draw=black!52,line width=.38pt] (v4) -- (v13);
\draw[draw=black!52,line width=.38pt] (v5) -- (v9);
\draw[draw=witnessred,line width=.8pt] (v5) -- (v10);
\draw[draw=black!52,line width=.38pt] (v5) -- (v11);
\draw[draw=witnessred,line width=.8pt] (v5) -- (v12);
\draw[draw=black!52,line width=.38pt] (v6) -- (v7);
\draw[draw=black!52,line width=.38pt] (v6) -- (v10);
\draw[draw=black!52,line width=.38pt] (v6) -- (v11);
\draw[draw=black!52,line width=.38pt] (v6) -- (v12);
\draw[draw=black!52,line width=.38pt] (v7) -- (v8);
\draw[draw=black!52,line width=.38pt] (v7) -- (v11);
\draw[draw=black!52,line width=.38pt] (v7) -- (v13);
\draw[draw=black!52,line width=.38pt] (v8) -- (v11);
\draw[draw=witnessred,line width=.8pt] (v8) -- (v12);
\draw[draw=black!52,line width=.38pt] (v9) -- (v11);
\draw[draw=black!52,line width=.38pt] (v9) -- (v12);
\draw[draw=black!52,line width=.38pt] (v9) -- (v13);
\draw[draw=black!52,line width=.38pt] (v10) -- (v11);
\draw[draw=witnessred,line width=.8pt] (v10) -- (v13);
\node[jfef vertex,draw=witnessblue,line width=.8pt] at (v0) {$0$};
\node[jfef vertex,draw=witnessblue,line width=.8pt] at (v1) {$1$};
\node[jfef vertex] at (v2) {$2$};
\node[jfef vertex] at (v3) {$3$};
\node[jfef vertex] at (v4) {$4$};
\node[jfef vertex,fill=orange!12] at (v5) {$5$};
\node[jfef vertex] at (v6) {$6$};
\node[jfef vertex] at (v7) {$7$};
\node[jfef vertex,fill=orange!12] at (v8) {$8$};
\node[jfef vertex] at (v9) {$9$};
\node[jfef vertex,fill=orange!12] at (v10) {$10$};
\node[jfef vertex] at (v11) {$11$};
\node[jfef vertex,fill=orange!12] at (v12) {$12$};
\node[jfef vertex,fill=orange!12] at (v13) {$13$};
\end{tikzpicture}
\caption{The solid edges form $\Ga$. With the dashed edge $\{0,1\}$,
the resulting graph belongs to $\IS$. The thick solid edges trace the induced path
$8$--$12$--$5$--$10$--$13$, the common neighbourhood of $0$ and $1$.
Thus the dashed edge is admissible for deletion.}\label{fig:G14}
\end{figure}
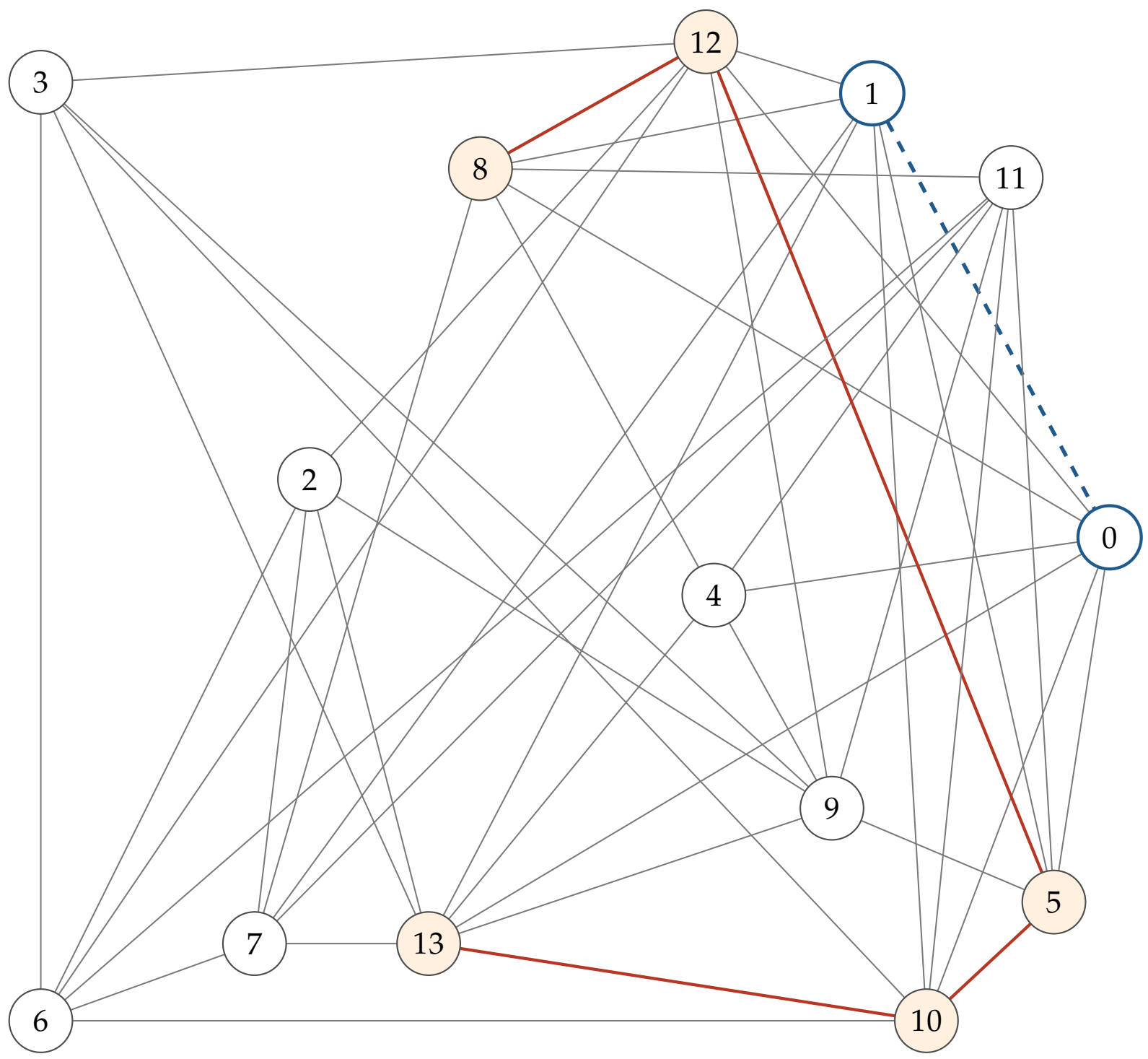

Two other graphs in Table~\ref{tab:witness-adjacency} admit the same
kind of extension, with a different induced path in each case.

\begin{proposition}\label{prop:other-seeds}
The graphs $\Gb$ and $\Gc$ are $K_4$-free members of
$\I\setminus\IS$, each with $14$ vertices and $45$ edges.
The extensions $\Gb+\{0,3\}$ and $\Gc+\{2,3\}$ belong to $\IS$,
and the added edges are admissible. The graphs $\Ga,\Gb,\Gc$ are
pairwise non-isomorphic.
\end{proposition}

\begin{proof}
In $\Gb$ the common neighbourhood of $0$ and $3$ is the induced
path $7\text{--}13\text{--}9\text{--}12$. In $\Gc$ the corresponding
path for $2$ and $3$ is $6\text{--}13\text{--}7\text{--}10\text{--}11$.
Thus both additions are admissible. Each extension reduces to vertex
$13$ by the deletions in Table~\ref{tab:variant-reductions}. All the
neighbourhoods encountered there are nonempty trees, so the extensions
belong to $\IS$. Deleting the added edge is an I-transformation in
each case, which gives $\Gb,\Gc\in\I$.

As in Theorem~\ref{thm:B}, neither original graph admits an admissible
deletion. Every vertex neighbourhood has positive first Betti number
and every edge neighbourhood is disconnected
(Table~\ref{tab:witness-links}). The same table distinguishes the two
graphs. In $\Gb$ there is an edge whose endpoints have four common neighbours,
whereas $\Gc$ has no such edge. Their edge counts distinguish both from
$\Ga$.
\end{proof}

The common neighbours of the endpoints of an edge correspond
bijectively to the triangles containing that edge. Thus the local
obstructions to a graph reduction also constrain simplicial collapses.

\begin{proposition}\label{prop:noncollapsible}
The clique complexes of $\Ga,\Gb,\Gc$ are contractible and
non-collapsible.
\end{proposition}

\begin{proof}
Contractibility follows from membership in $\I$. Each graph contains a triangle and is $K_4$-free,
so its clique complex is two-dimensional. Every
edge belongs to at least two triangles by Table~\ref{tab:witness-links},
and every vertex lies on an edge. Thus each vertex, too, lies in at
least two triangles. None of the three complexes has a free face, so no elementary
collapse can begin.
\end{proof}

\section{Algorithms and the enumeration of graph classes}\label{sec:algorithms}

The examples of Section~\ref{sec:examples} distinguish
$\mathbb F_2$-acyclicity, I-contractibility and membership in $\IS$.
To determine the first order at which these distinctions occur, we
must also account for every graph of smaller order. The homology of
the clique complex and the possible first deletions provide two ways
to do so. Homology detects obstructions to I-contractibility, while
the deletion recurrences decide membership in $\IS$, $\ISV$ and the
dismantling classes. For graphs outside $\IS$, the one-edge
constructions of Section~\ref{sec:fourteen} suggest a further test
for membership in $\I$.

\subsection{Cliques and homology}

Fix a linear order on the vertices of a nonempty graph $G$. Starting
with each singleton, extend the current clique in all possible ways
by adjoining a common neighbour larger than its largest vertex.
Each clique occurs exactly once, since its vertices in increasing order
determine the successive extensions. The cliques of size $i+1$ thus
obtained are the $i$-simplices of $\cl(G)$.

Let $C_i$ denote the integral simplicial chain group of $\cl(G)$,
with each clique oriented by increasing vertex order. For $i\ge 1$ the boundary map is
\[
 \partial_i[v_0,\ldots,v_i]
   =\sum_{j=0}^i(-1)^j[v_0,\ldots,\widehat v_j,\ldots,v_i].
\]
Adjoin the augmentation $\partial_0:C_0\to\ZZ$, which sends every
vertex to $1$. Write $f_i=\operatorname{rank} C_i$ and let $r_i(R)$
denote the rank of $\partial_i$ over the field $R$. Since $\operatorname{im}\partial_{i+1}\subseteq\ker\partial_i$,
rank--nullity gives
\begin{equation}\label{eq:betti}
 \dim_R\widetilde H_i(G;R)=f_i-r_i(R)-r_{i+1}(R)
 \qquad(i\ge0),
\end{equation}
where $r_0(R)=1$ and $r_i(R)=0$ for $i>\dim\cl(G)$.
A graph in $\mathcal A_{\mathbb F_2}$ has Euler characteristic one.
For graphs of this Euler characteristic, row reduction over
$\mathbb F_2$ gives the ranks in~\eqref{eq:betti} and hence decides
membership.

The distinction between acyclicity over $\mathbb F_2$ and over $\ZZ$
depends on the integral boundary maps. Let $d_{i+1,j}$
denote the nonzero invariant factors in the Smith normal form of
$\partial_{i+1}$. Then
\begin{equation}\label{eq:integral-homology}
 \widetilde H_i(G;\ZZ)\cong
 \ZZ^{\,f_i-r_i(\mathbb Q)-r_{i+1}(\mathbb Q)}
 \oplus\bigoplus_{d_{i+1,j}>1}\ZZ/d_{i+1,j}.
\end{equation}
Indeed, the augmented chain complex gives a short exact sequence
\[
 0\longrightarrow\widetilde H_i(G;\ZZ)
 \longrightarrow C_i/\operatorname{im}\partial_{i+1}
 \longrightarrow\operatorname{im}\partial_i\longrightarrow0.
\]
The last group is free abelian, so the sequence splits. The torsion subgroup of
the middle group is therefore isomorphic to that of
$\widetilde H_i(G;\ZZ)$, and the ranks of the two boundary maps give
the free rank in~\eqref{eq:integral-homology}. Applying this formula
in each degree decides membership in $\mathcal A_{\ZZ}$.

\subsection{Recursive criteria for graph reductions}

A reduction to $K_1$ records more than the homology of $\cl(G)$.
Its first step singles out both a neighbourhood and a smaller graph
that themselves admit reductions. To express this recursion, let $S$
be the Boolean function on graphs defined by
$S(\emptyset)=\mathrm{false}$, $S(K_1)=\mathrm{true}$ and, for
$|V(G)|\ge 2$, by the recurrence
\begin{align}
 S(G)={}&\bigvee_{v\in V(G)}
 \bigl(S(\Nb_G(v))\land S(G-v)\bigr)\notag\\
 &{}\lor\bigvee_{uv\in E(G)}
 \bigl(S(\Nb_G(u,v))\land S(G-uv)\bigr).
 \label{eq:strong-rec}
\end{align}
Each neighbourhood in~\eqref{eq:strong-rec} is taken in the graph
at the current step of the reduction.

\begin{proposition}\label{prop:IS-correct}
The recurrence~\eqref{eq:strong-rec} terminates on every graph, and
$S(G)$ is true if and only if $G\in\IS$. Whenever $S(G)$ is true,
the recurrence yields a reduction of $G$ to $K_1$ and reductions of
all the neighbourhoods used to justify the deletions.
\end{proposition}
\begin{proof}
Each term involves only graphs with smaller $|V|+|E|$, and there
are finitely many terms. Thus the recursion terminates. By induction
on $|V(G)|+|E(G)|$, a given term is true
exactly when the neighbourhood and the remaining graph both belong
to $\IS$. Definition~\ref{def:IS} gives the asserted equivalence.
Choosing one true term in each nontrivial recursive call gives the
required reductions.
\end{proof}

For reductions that delete only vertices, the remaining graph always
has the form $G[U]$ for some $U\subseteq V(G)$. If $v\in U$, write
$N_U(v)=U\cap V(\Nb_G(v))$ for the set of its neighbours in $U$.
We define Boolean functions $W$ and $D_k$ on the subsets of $V(G)$
for membership in $\ISV$ and $k$-Dis, respectively. Each function
is false on the empty set and true on singletons. For nonempty $U$,
$D_{-1}(U)$ is true exactly when $G[U]$ is a cone. For $|U|\ge2$,
the recurrences are
\begin{align}
 W(U)&=\bigvee_{v\in U}
     \bigl(W(N_U(v))\land W(U\setminus\{v\})\bigr),
     \label{eq:ivs-rec}\\
 D_k(U)&=\bigvee_{v\in U}
     \bigl(D_{k-1}(N_U(v))\land D_k(U\setminus\{v\})\bigr),
     \qquad k\ge0.\label{eq:kdis-rec}
\end{align}

\begin{proposition}\label{prop:vertex-correct}
The recurrences~\eqref{eq:ivs-rec} and~\eqref{eq:kdis-rec} decide
whether $G[U]$ belongs to $\ISV$ and to $k$-Dis, respectively. For a
fixed upper bound on $k$, both recurrences can be evaluated in
increasing order of $|U|$.
\end{proposition}
\begin{proof}
The initial values agree with Definitions~\ref{def:IS}
and~\ref{def:kdis}. For $|U|\ge2$, the two sets in each term are
proper subsets of $U$. By induction on $|U|$, a vertex $v$ contributes
a true term to~\eqref{eq:ivs-rec} precisely when its neighbourhood and
the graph remaining after its deletion both belong to $\ISV$.
For~\eqref{eq:kdis-rec} these two graphs must belong to $(k-1)$-Dis
and $k$-Dis, respectively. In each case the disjunction runs through
all possible first vertices.
\end{proof}

To decide $0$-dismantlability one may delete dominated vertices in
any order \cite[Section~4.1]{FieuxJouve2020}. Indeed, a vertex
$w\ne v$ dominates $v$ if $N_G[v]\subseteq N_G[w]$, equivalently
if $\Nb_G(v)$ is a cone with apex $w$. Suppose $u$ and $v$ are
distinct dominated vertices. If they dominate each other, their
transposition is an automorphism, so $G-u$ and $G-v$ are isomorphic.
Otherwise $v$ remains dominated in $G-u$. Any dominator of $v$
distinct from $u$ still dominates it there. If $u$ dominates $v$,
a dominator of $u$ is distinct from $v$ and also dominates $v$ by
transitivity. The same argument with $u$ and $v$ interchanged shows
that either deletion can be followed by the other, giving
$G-\{u,v\}$ in both orders. Induction on $|V(G)|$ now shows that
whenever one maximal sequence of dominated-vertex deletions ends at
$K_1$, every such sequence does. We use this criterion to count
$0$-dismantlable graphs. The choice of deletion matters for $\IS$.
Each of the extensions in Section~\ref{sec:fourteen} belongs to
$\IS$, yet admits an edge deletion that leaves a graph outside $\IS$.

\subsection{A sufficient condition for I-contractibility}

A graph outside $\IS$ may acquire a reduction to $K_1$ after a single
admissible edge addition, as $\Ga$ does. More generally, suppose that
\begin{equation}\label{eq:I-path}
 G=G_0,\qquad G_j=G_{j-1}+u_jv_j\quad(1\le j\le t),
 \qquad G_t\in\IS,
\end{equation}
and $\Nb_{G_{j-1}}(u_j,v_j)\in\IS$ for every $j$. Then $G\in\I$.
Indeed, $G_t\in\IS\subseteq\I$, and the addition of $u_jv_j$ leaves
the common neighbourhood of $u_j$ and $v_j$ unchanged. Deleting the
added edges in reverse order therefore gives a sequence of
I-transformations from $G_t$ to $G$.

For a fixed bound $t$, we examine the graphs reached from $G$ by at
most $t$ admissible edge additions, applying
Proposition~\ref{prop:IS-correct} to each. A graph in $\IS$ at the
end of such a sequence certifies $G\in\I$. If there is no such
sequence within the bound, membership in $\I$ may still be settled
by other means. In particular, nonzero reduced integral homology
excludes $G$ from $\I$ by~\eqref{eq:homology-inclusions}.

\subsection{Representatives up to isomorphism}\label{sec:representatives}

To apply the tests to all graphs of a prescribed order, it suffices
to generate one graph from each isomorphism class. McKay's method of
canonical augmentation does this by choosing which vertex may be added
last \cite[Section~4]{McKay1998}. The following construction illustrates
the principle used by \texttt{geng}, the generator employed in
Section~\ref{sec:searches}. For each graph $H$ of order at least
two, choose a vertex orbit $O(H)$ under $\operatorname{Aut}(H)$ so
that each isomorphism $H\to H'$ carries $O(H)$ to $O(H')$. For example, take the orbit of the vertex assigned label $1$ in a canonical labelling with labels $1,\ldots,|V(H)|$. Starting
from $K_1$, form augmentations of each retained graph $G$ by adjoining
a vertex $x$. For its neighbourhood choose one representative from each
orbit of subsets of $V(G)$ under $\operatorname{Aut}(G)$. Retain an
augmentation $H$ precisely when $x\in O(H)$. In this construction,
connectedness is tested only at the prescribed order, so intermediate
graphs may be disconnected. To enumerate $K_4$-free graphs, we also discard every
augmentation containing a $K_4$.

\begin{proposition}\label{prop:enumeration-correct}
Canonical vertex augmentation followed by the connectedness test gives
one representative of each isomorphism class of connected graphs of
the prescribed order. With the additional restriction on augmentations,
it gives precisely the representatives of the connected $K_4$-free
graphs of that order.
\end{proposition}
\begin{proof}
Before imposing connectedness, we argue by induction on the order.
The unique graph of order one is $K_1$. Given $H$ of larger order,
choose $x\in O(H)$. By induction, the construction contains a
representative of $H-x$. Augmenting that representative with the
neighbourhood of $x$, up to the action of its automorphism group,
gives a copy of $H$. Since the chosen orbit is invariant under
isomorphism, this augmentation is retained.

Suppose two retained augmentations are isomorphic. Their new vertices
belong to the respective chosen orbits, so we can compose the
isomorphism with an automorphism of the target to match the two new
vertices. The restriction to the remaining vertices identifies the
graphs before augmentation. By induction these have the same
representative, and the two neighbourhoods belong to the same orbit
of subsets of its vertex set. Only one such augmentation was made.
This proves uniqueness. At the prescribed order the connectedness
test retains precisely the connected representatives. Since induced
subgraphs of a $K_4$-free graph are $K_4$-free, the additional
restriction retains every predecessor needed in the same argument.
\end{proof}

Applying these tests to one representative of each class determines
the counts for $\mathcal A_{\mathbb F_2}$, $\IS$, $\ISV$ and the
dismantling levels. By~\eqref{eq:homology-inclusions}, only the
representatives in $\mathcal A_{\mathbb F_2}\setminus\IS$ require
further analysis to decide membership in $\I$.

\subsection{Reduction of the order-fourteen classification}\label{sec:reduction}

At the minimum order of a graph in $\I\setminus\IS$, admissible
deletions put a further constraint on the classification. Each such
deletion preserves I-contractibility. It also leaves the graph outside
$\IS$, since a reduction of the smaller graph to $K_1$ would extend
to one of the original graph. Thus a vertex deletion would contradict
minimality, while successive edge deletions eventually reach a graph
with no admissible deletion. We shall also need the corresponding
condition on either side of a \defi{cut vertex}, a vertex whose deletion increases the number of connected components.

For a graph $G$ of order at least two, consider the following conditions.
\begin{description}[font=\normalfont,leftmargin=3.7em]
\item[\hypertarget{cond:no-deletion}{\textup{(ND)}}] $G$ is $\mathbb F_2$-acyclic and has no admissible vertex or edge.
\item[\hypertarget{cond:no-deletion-away}{\ensuremath{(\mathrm{ND}_u)}}] For a specified vertex $u$, the graph $G$ is $\mathbb F_2$-acyclic
and has no admissible edge and no admissible vertex other than possibly $u$.
\end{description}

\begin{lemma}\label{lem:stuck}
Let $\mathcal H$ be closed under vertex and edge deletion, and let $n$
be the smallest order attained in $\mathcal H\cap(\I\setminus\IS)$.
Every graph of order $n$ in this intersection admits a sequence of
admissible edge deletions ending at a graph in the intersection with
no admissible vertex or edge.
\end{lemma}
\begin{proof}
An admissible vertex deletion would give a graph of smaller order in
$\mathcal H\cap(\I\setminus\IS)$, contrary to the choice of $n$.
Each admissible edge deletion preserves both the order and membership
in this intersection. Since each step removes an edge, a maximal
sequence of such deletions ends at a graph with no admissible edge.
It has no admissible vertex by the same minimality argument and is
$\mathbb F_2$-acyclic because it belongs to $\I$.
\end{proof}

If $u$ is a cut vertex of a graph satisfying \hyperlink{cond:no-deletion}{\textup{(ND)}}, then a smaller induced
subgraph containing $u$ satisfies \hyperlink{cond:no-deletion-away}{\ensuremath{(\mathrm{ND}_u)}}.

\begin{lemma}\label{lem:2conn}
Suppose $G$ satisfies condition~\hyperlink{cond:no-deletion}{\textup{(ND)}} and has a cut vertex $u$. Then $G$ is the
union of two induced subgraphs $A$ and $B$, each of order at least
three, meeting only at $u$. Both satisfy \hyperlink{cond:no-deletion-away}{\ensuremath{(\mathrm{ND}_u)}}. Consequently, if no graph $H$ of order at most $m$ in a hereditary
class satisfies \hyperlink{cond:no-deletion-away}{\ensuremath{(\mathrm{ND}_w)}} for any $w\in V(H)$, then every graph in the
class of order at most $2m$ satisfying \hyperlink{cond:no-deletion}{\textup{(ND)}} is $2$-connected.
\end{lemma}
\begin{proof}
Partition the components of $G-u$ into two nonempty collections. Let
$A$ and $B$ be the subgraphs induced by their respective vertex sets
together with $u$. Then $\cl(G)=\cl(A)\cup\cl(B)$ and
$\cl(A)\cap\cl(B)=\{u\}$. The reduced Mayer--Vietoris sequence gives
\[
 \widetilde H_i(G;\mathbb F_2)\cong
 \widetilde H_i(A;\mathbb F_2)\oplus
 \widetilde H_i(B;\mathbb F_2).
\]
Thus $A$ and $B$ are $\mathbb F_2$-acyclic. Since there are no edges
between $A-u$ and $B-u$, every vertex other than $u$ has the same
neighbourhood in $G$ as in its induced subgraph. The common
neighbourhood of the endpoints of any edge is likewise unchanged,
since at least one endpoint differs from $u$. This proves
\hyperlink{cond:no-deletion-away}{\ensuremath{(\mathrm{ND}_u)}}. If either subgraph had order two, its vertex distinct
from $u$ would be a leaf of $G$ and hence admissible. Thus both
subgraphs have order at least three. In a hereditary class they also
belong to the class, so the hypothesis in the second assertion forces
both orders to exceed $m$. Their union would then have order at least
$2(m+1)-1=2m+1$. This rules out a cut vertex in the stated range.
A graph satisfying \hyperlink{cond:no-deletion}{\textup{(ND)}} is connected by acyclicity and cannot
have order two, since either vertex of $K_2$ is admissible.
\end{proof}

Once a cut vertex has been ruled out, acyclicity and the absence of
admissible deletions impose bounds on the vertex degrees.

\begin{lemma}\label{lem:degrees}
Let $G$ be a $2$-connected, $\mathbb F_2$-acyclic graph with no
admissible vertex or edge. Every edge of $G$ lies in at least two
triangles. The minimum degree $\delta(G)$ and maximum degree $\Delta(G)$ of $G$ satisfy
\[
 4\le\delta(G)\le\Delta(G)\le |V(G)|-2.
\]
\end{lemma}
\begin{proof}
Suppose an edge belongs to no triangle. Since $G$ is $2$-connected,
the edge lies on a cycle.
Over $\mathbb F_2$, the coefficient of this edge in the boundary of
every $2$-chain is zero, whereas its coefficient in the cycle is one.
The cycle therefore represents a nonzero class in $H_1(G;\mathbb F_2)$,
contrary to acyclicity.
If an edge lies in exactly one triangle, the common neighbourhood
of its endpoints is $K_1$, so the edge is admissible. Hence every
edge lies in at least two triangles.
For $H=\Nb_G(v)$, every vertex $w$ of $H$ consequently has at least
two neighbours in $H$. If $|V(H)|\le3$, this forces $H=K_3$, making
$v$ admissible. Thus every degree is at least four. A universal
vertex would make $G$ a cone and hence a member of $\IS$, so no
degree is $|V(G)|-1$.
\end{proof}

In a $K_4$-free graph, the neighbourhood of a vertex contains no
triangle. This makes it possible to recognise the admissible
deletions directly from the neighbourhoods.
\begin{proposition}\label{rem:K4-local}
In a $K_4$-free graph, a vertex is $\IS$-admissible exactly when its
open neighbourhood is a nonempty tree. An edge is $\IS$-admissible
if and only if it belongs to exactly one triangle.
\end{proposition}
\begin{proof}
A vertex neighbourhood is triangle-free. A triangle-free graph in
$\IS$ has a contractible clique complex of dimension at most one,
so it is a tree. Conversely, a nonempty tree reduces by leaf deletions
to $K_1$. The common neighbourhood of the endpoints of an edge is an edgeless
graph. It belongs to $\IS$ exactly when it has one vertex. Its vertices
correspond to the triangles containing the edge.
\end{proof}

It remains to recover the graphs that admit edge deletions. Reversing
the sequences in Lemma~\ref{lem:stuck} suggests the following
construction. Start with a finite collection of
$\mathbb F_2$-acyclic graphs of order at least two with no admissible
vertex or edge. Each lies outside $\IS$, since a reduction to $K_1$
requires a first deletion. From each graph $H$, form all graphs $H+uv$
where $uv$ is a missing edge and $\Nb_H(u,v)\in\IS$. Retain those
outside $\IS$ and repeat the construction from each newly obtained
isomorphism class.

\begin{proposition}\label{prop:closure-correct}
This construction terminates. Including the initial graphs, it gives
exactly the isomorphism classes outside $\IS$ that can be obtained by
admissible edge additions from the initial collection. The same
conclusion holds within the $K_4$-free graphs if graphs containing
a $K_4$ are discarded throughout.
\end{proposition}
\begin{proof}
Each sequence of additions has a fixed vertex set and strictly
increasing edge count. Hence a finite initial collection gives only
finitely many graphs. Every graph retained belongs to a sequence
of the stated kind. For the converse, observe that an admissible
addition to $H\in\IS$ gives $H+uv\in\IS$, since deleting $uv$ and
then reducing $H$ gives a reduction to $K_1$. Thus every graph in a
sequence ending outside $\IS$ also lies outside $\IS$. By induction
on the length of the sequence, the construction retains a
representative at every step. Indeed, an isomorphism preserves both
nonadjacency and the common neighbourhood up to isomorphism, so all
admissible additions occur from any representative. Finally, a $K_4$
persists under edge addition. Every sequence ending at a $K_4$-free
graph therefore consists entirely of $K_4$-free graphs.
\end{proof}

\section{Cardinalities and minimum orders}\label{sec:searches}\label{sec:table}

To prove the minimality of the examples in Section~\ref{sec:examples},
we use the membership tests of Section~\ref{sec:algorithms} to count
the graphs in each family. Table~\ref{tab:card} covers all connected graphs of
orders $3$ through $13$ and all connected $K_4$-free graphs of order
$14$. Every graph in $\mathcal A_{\mathbb F_2}$ is connected, so
connectedness imposes no further restriction on the acyclic graphs.

We obtained one representative per isomorphism class with McKay's
generator \texttt{geng} \cite[Section~4]{McKay1998} from the nauty
package \cite{McKayPiperno2014}. We used nauty 2.8.9 through order
$11$ and 2.7r4 for orders $12$ and $13$. At order $14$ the $K_4$-free option of nauty 2.9.3 served for
two enumerations. The first covers all connected $K_4$-free graphs
and supplies the last row of the table. In the second we used the
options for $2$-connectedness and the degree bounds of Lemma~\ref{lem:degrees} for
the proof of Theorem~\ref{thm:classification14}. For orders $12$ through
$14$, the counts come from the same family of classification programs.
The repository~\cite{GraphData2026} records the generation parameters and aggregate counts
for the partitions of each enumeration. Through order $9$, two
implementations agreed on the classification of each graph. For the
larger enumerations, separate programs check the
homology and reductions of the exceptional graphs in the catalogues.

In developing the enumeration software we used Claude (Anthropic) to
implement the C classifier, the programs for running and recording
the enumerations, and the separate routines for checking homology and
replaying reduction certificates. We used Codex (OpenAI) for a
separate enumeration at small orders. The repository~\cite{GraphData2026} also includes the
script that draws the labelled graphs and the unfolded disc in
Figures~\ref{fig:G13} and~\ref{fig:G14} from the recorded incidences
and vertex coordinates.

For the larger enumerations we first tested the Euler characteristic,
computing the boundary ranks over $\mathbb F_2$ only when $\chi=1$.
Successive deletions of dominated vertices decide $0$-dismantlability.
A reduction of this kind to $K_1$ also puts the graph in both $\ISV$
and $\IS$. For the remaining $\mathbb F_2$-acyclic graphs, we tested
membership first in $\ISV$ and then, if necessary, in $\IS$.
The graphs left in $\mathcal A_{\mathbb F_2}\setminus\IS$ determine
the entries for $\I$ and $\mathcal A_{\ZZ}$, as we show below.

Two equalities account for the combined columns.
Proposition~\ref{prop:table-integral} identifies $\mathcal A_{\ZZ}$
with $\I$ in the range of each row. Through order $11$ the equal counts for $\Dis{1}$ and $\ISV$,
together with~\eqref{eq:dismantling-inclusions}, determine the entries
for $\Dis{2}$ and $\Dis{3}$. At orders $12$ and $13$ the counts for
$\Dis{2}$, $\Dis{3}$ and $\ISV$ agree. For $K_4$-free graphs the
corresponding equality follows from the tree-neighbourhood criterion
in Proposition~\ref{rem:K4-local}, as shown below. The cone counts
follow by adjoining a universal vertex to a graph of one smaller order.

\begin{table}[!t]
\centering
\caption{Numbers of isomorphism classes of connected graphs.
Equalities in the column headings hold under the restrictions of each row.
The cones are the $(-1)$-dismantlable graphs.}\label{tab:card}
\begingroup
\footnotesize
\setlength{\tabcolsep}{1pt}
\renewcommand{\arraystretch}{1.07}
\sisetup{group-separator={,}}
\begin{tabular*}{\linewidth}{@{\extracolsep{\fill}}rrrrr@{}}
\toprule
\multicolumn{5}{c}{All connected graphs for $3\le n\le13$ and connected $K_4$-free graphs for $n=14$ ($\dagger$)} \\
\midrule
$n$ & Connected & $\chi=1$ & $\mathcal A_{\mathbb F_2}$ & $\mathcal A_{\ZZ}=\I$ \\
\midrule
3 & \num{2} & \num{2} & \num{2} & \num{2} \\
4 & \num{6} & \num{5} & \num{5} & \num{5} \\
5 & \num{21} & \num{16} & \num{16} & \num{16} \\
6 & \num{112} & \num{68} & \num{68} & \num{68} \\
7 & \num{853} & \num{404} & \num{403} & \num{403} \\
8 & \num{11117} & \num{3837} & \num{3798} & \num{3798} \\
9 & \num{261080} & \num{68458} & \num{66167} & \num{66167} \\
10 & \num{11716571} & \num{2571263} & \num{2343148} & \num{2343148} \\
11 & \num{1006700565} & \num{204131202} & \num{166978401} & \num{166978401} \\
12 & \num{164059830476} & \num{32482574568} & \num{22666537431} & \num{22666537431} \\
13 & \num{50335907869219} & \num{9859955384222} & \num{5619022190286} & \num{5619022190283} \\
\addlinespace[2pt]\cmidrule(lr){1-5}
14$^{\dagger}$ & \num{309317860269934} & \num{6112638767338} & \num{755877673752} & \num{755877673127} \\
\bottomrule
\end{tabular*}
\par\smallskip
\begin{tabular*}{\linewidth}{@{\extracolsep{\fill}}rrrrrr@{}}
\toprule
$n$ & $\IS$ & $2\text{-Dis}=3\text{-Dis}=\ISV$ & $1$-Dis & $0$-Dis & Cones \\
\cmidrule(r){1-1}\cmidrule(lr){2-2}\cmidrule(lr){3-3}\cmidrule(lr){4-4}\cmidrule(lr){5-5}\cmidrule(l){6-6}
3 & \num{2} & \num{2} & \num{2} & \num{2} & \num{2} \\
4 & \num{5} & \num{5} & \num{5} & \num{5} & \num{4} \\
5 & \num{16} & \num{16} & \num{16} & \num{16} & \num{11} \\
6 & \num{68} & \num{68} & \num{68} & \num{68} & \num{34} \\
7 & \num{403} & \num{403} & \num{403} & \num{403} & \num{156} \\
8 & \num{3798} & \num{3798} & \num{3798} & \num{3791} & \num{1044} \\
9 & \num{66167} & \num{66167} & \num{66167} & \num{65561} & \num{12346} \\
10 & \num{2343148} & \num{2343148} & \num{2343148} & \num{2258314} & \num{274668} \\
11 & \num{166978401} & \num{166978389} & \num{166978389} & \num{151084672} & \num{12005168} \\
12 & \num{22666537431} & \num{22666525713} & \num{22666525693} & \num{18510237804} & \num{1018997864} \\
13 & \num{5619022190283} & \num{5619007241580} & \num{5619007205823} & \num{4015494549880} & \num{165091172592} \\
\addlinespace[2pt]\cmidrule(lr){1-6}
14$^{\dagger}$ & \num{755877673119} & \num{755266302709} & \num{755266302709} & \num{81156275851} & \num{20797002} \\
\bottomrule
\end{tabular*}
\endgroup

\end{table}

\subsection{\texorpdfstring{$\mathbb F_2$}{F2}-acyclic graphs of order at most thirteen}

Through order $12$, acyclicity over $\mathbb F_2$ already implies the
existence of a reduction to $K_1$. The three graphs of
Theorem~\ref{thm:A} account for the first exceptions.

\begin{theorem}\label{thm:classification13}
Up to isomorphism, the only $\mathbb F_2$-acyclic graphs on at most
$13$ vertices outside $\IS$ are $\Gt,\Gt',\Gt''$. Every graph in
$\I$ on at most $13$ vertices belongs to $\IS$. Among the
$\mathbb F_2$-acyclic graphs of order $13$, $\Gt$ is the only one
up to isomorphism with no admissible vertex or edge.
\end{theorem}
\begin{proof}
At orders $3$ through
$12$, the inclusion $\IS\subseteq\mathcal A_{\mathbb F_2}$ and the
equal counts in Table~\ref{tab:card} give equality of the two families.
The only connected graphs of orders one and two are $K_1$ and $K_2$,
both of which belong to $\IS$. At order $13$ the counts differ by
\[
 \RowThirteenAcyclic-\RowThirteenIS=3.
\]
By Theorem~\ref{thm:A}, the three graphs $\Gt,\Gt',\Gt''$ are
$\mathbb F_2$-acyclic, and their nonzero first integral homology
excludes them from $\I$ and hence from $\IS$. They have different
numbers of edges and therefore account for all three isomorphism
classes in the difference. Since $\I\subseteq\mathcal A_{\mathbb F_2}$
and none of the three belongs to $\I$, every graph in $\I$ of order
at most $13$ belongs to $\IS$.

An $\mathbb F_2$-acyclic graph of order $13$ with no admissible vertex
or edge lies outside $\IS$, since a reduction to $K_1$ would begin
with such a deletion. Table~\ref{tab:witness-links} shows that $\Gt$
has this property. In each of $\Gt'$ and $\Gt''$ the last edge added
in the construction is admissible for deletion.
\end{proof}

\subsection{The \texorpdfstring{$K_4$}{K4}-free graphs of order fourteen}

Theorem~\ref{thm:classification13} and the example $\Ga$ place the
first graphs in $\I\setminus\IS$ at order $14$. For those without
a $K_4$, Lemma~\ref{lem:stuck} reduces the classification to a smaller
family. We first find the graphs with no admissible vertex or edge,
then recover the others by reversing edge deletions.

\begin{theorem}\label{thm:classification14}
Up to isomorphism, there are $51$ $\mathbb F_2$-acyclic $K_4$-free
graphs of order $14$ with no admissible vertex or edge. Exactly three
of them belong to $\I$, namely $\Ga,\Gb,\Gc$. Each of the other $48$
has first integral homology $\ZZ/3$ and vanishing higher integral
homology. There are exactly eight isomorphism classes of $K_4$-free
graphs of order $14$ in $\I\setminus\IS$. Among these graphs the
minimum number of edges is $44$, attained by $\Ga$.
\end{theorem}
\begin{proof}
Let $G$ be an $\mathbb F_2$-acyclic $K_4$-free graph of order $14$
with no admissible vertex or edge. Suppose $G$ has a cut vertex $u$.
By Lemma~\ref{lem:2conn}, $G$ is the union of two induced subgraphs
meeting at $u$ whose orders sum to $15$. Each is $\mathbb F_2$-acyclic,
has no admissible edge and has no admissible vertex other than
possibly $u$. At least one of the two subgraphs has at most seven
vertices. At orders two through seven, there are respectively
$1,2,5,17,82,536$ isomorphism classes of connected $K_4$-free graphs.
We tested each representative for $\mathbb F_2$-acyclicity and
applied the neighbourhood criteria of Proposition~\ref{rem:K4-local}.
For every choice of $u$ in each $\mathbb F_2$-acyclic representative, the tests identify
an admissible edge or an admissible vertex distinct from $u$.
Thus $G$ is $2$-connected and, by Lemma~\ref{lem:degrees},
$4\le\delta(G)\le\Delta(G)\le12$.

We thus enumerated the $2$-connected $K_4$-free graphs of order $14$
with minimum degree at least four and maximum degree at most twelve.
There are $\num{64856470381820}$ such graphs up to isomorphism. Among
the graphs of Euler characteristic one, $\num{13536597}$ had no
admissible vertex or edge, as determined by
Proposition~\ref{rem:K4-local}. Each clique complex has dimension at
most two, so it remained to compute the ranks of $\partial_1$ and
$\partial_2$ over $\mathbb F_2$. By~\eqref{eq:betti}, exactly $51$ of
the $\num{13536597}$ representatives had all their reduced Betti
numbers equal to zero. Table~1 of the supplementary material lists these $51$ graphs.

Theorem~\ref{thm:B} and Proposition~\ref{prop:other-seeds} show that
$\Ga,\Gb,\Gc$ belong to $\I$. For each of the remaining $48$ graphs,
formula~\eqref{eq:integral-homology} gives free rank zero for the first
integral homology group, and the only nonunit Smith factor of
$\partial_2$ is $3$. Thus $H_1\cong\ZZ/3$, which excludes these graphs
from $\I$. Their clique complexes are connected, have dimension at most
two and have Euler characteristic one. Hence $H_2$ is free of rank zero
and all higher homology vanishes.

Since $14$ is the minimum order for a graph in $\I\setminus\IS$,
Lemma~\ref{lem:stuck} now gives the
rest of the classification. Every $K_4$-free graph of this order in
$\I\setminus\IS$ can be transformed into one of $\Ga,\Gb,\Gc$ by
admissible edge deletions. To reverse these deletions, consider a missing edge $uv$. Its
addition is admissible and preserves $K_4$-freeness precisely when
$u$ and $v$ have a unique common neighbour. Indeed, to avoid a $K_4$
their common neighbourhood must be edgeless, and the only edgeless
graph in $\IS$ is $K_1$.

From Table~\ref{tab:witness-adjacency}, the only such pairs in $\Ga$
are $\{3,8\}$ and $\{4,6\}$, with common neighbours $12$ and $11$,
respectively. In $\Gb$ the pairs are $\{0,6\}$ and $\{1,6\}$,
again with common neighbours $12$ and $11$, respectively. There
are none in $\Gc$. Every other non-adjacent pair has at least two
common neighbours, which remain common neighbours after any edge
addition. In each of $\Ga$ and $\Gb$, adding either of the two listed
edges leaves the common neighbourhood of the other unchanged.
Thus no other admissible addition can occur while the graph remains
$K_4$-free.

In every resulting graph, each original vertex neighbourhood retains
a cycle and each original edge lies in at least two triangles.
The added edges are precisely the edges in a single triangle.
By Proposition~\ref{rem:K4-local}, the only admissible deletions
therefore remove added edges. Every maximal sequence of admissible
deletions returns to the initial graph. Thus each extension lies
outside $\IS$ and extensions of non-isomorphic initial graphs are
non-isomorphic.

For $\Ga$, the two single-edge extensions are distinguished by the
degree of the third vertex in the unique triangle containing the
added edge. This degree is eight for vertex $12$ and seven for
vertex $11$. For $\Gb$, the automorphism
$(0\,1)(3\,4)(7\,8)(11\,12)$ interchanges $\{0,6\}$ and $\{1,6\}$,
so the single-edge extensions are isomorphic. Including each
initial graph, we obtain the following numbers of classes by edge count.
\[
\begin{array}{c|rrrr|r}
 \text{initial graph} & 44 & 45 & 46 & 47 & \text{total}\\ \hline
 \Ga & 1 & 2 & 1 & 0 & 4\\
 \Gb & 0 & 1 & 1 & 1 & 3\\
 \Gc & 0 & 1 & 0 & 0 & 1
\end{array}
\]
The three rows are disjoint by Proposition~\ref{prop:other-seeds}
and the uniqueness of the initial graph under deletion. Each of the eight graphs belongs to $\I$, since every admissible
edge addition preserves I-contractibility. By
Proposition~\ref{prop:closure-correct}, every $K_4$-free graph of
order $14$ in $\I\setminus\IS$ occurs in one of these rows. The
least edge count is $44$, attained by $\Ga$.
\end{proof}

For each of the eight graphs, the repository~\cite{GraphData2026} also gives an admissible
one-edge extension in $\IS$ and a reduction of that extension to $K_1$.

\subsection{Integral acyclicity and I-contractibility}

The same classification also determines integral acyclicity. In the
chain of inclusions
\begin{equation}\label{eq:table-chain}
 \IS\subseteq\I\subseteq\mathcal A_{\ZZ}
       \subseteq\mathcal A_{\mathbb F_2},
\end{equation}
only graphs in $\mathcal A_{\mathbb F_2}\setminus\IS$ can distinguish
$\I$ from $\mathcal A_{\ZZ}$. For these graphs the integral homology
is preserved by the admissible deletions used in the classification.

\begin{proposition}\label{prop:table-integral}
For connected graphs of order at most $13$ and connected $K_4$-free
graphs of order $14$, integral acyclicity is equivalent to
I-contractibility. In Table~\ref{tab:card}, the number of classes
with these equivalent properties equals the number of
$\mathbb F_2$-acyclic classes through order $12$. It is smaller by
three at order $13$ and by $625$ among the $K_4$-free graphs of order $14$.
\end{proposition}
\begin{proof}
For graphs on at most $13$ vertices, Theorem~\ref{thm:classification13} gives all the
graphs in $\mathcal A_{\mathbb F_2}\setminus\IS$. Each has nonzero
first integral homology by Theorem~\ref{thm:A}.
Now let $G$ be a $K_4$-free graph of order $14$ in
$\mathcal A_{\mathbb F_2}\setminus\IS$. Every admissible deletion preserves
the homotopy type of the clique complex. Furthermore, each graph obtained
from $G$ by admissible deletions lies outside $\IS$, since a reduction
of that graph to $K_1$ would extend to one of $G$.

Starting with $G$, delete admissible edges until either a vertex
becomes admissible for deletion or no admissible edge remains. If a
vertex is admissible at any stage, including in $G$ itself, delete it.
The resulting graph belongs to $\mathcal A_{\mathbb F_2}\setminus\IS$
and has order $13$, so it is isomorphic to one of the three graphs in
Theorem~\ref{thm:classification13}. Otherwise the edge deletions end at
a $K_4$-free spanning subgraph of $G$ with no admissible vertex or edge.
It is still $\mathbb F_2$-acyclic, so Theorem~\ref{thm:classification14}
identifies it either with one of the $48$ graphs with first integral
homology $\ZZ/3$ or with one of $\Ga,\Gb,\Gc$. In the latter case
the reverse sequence of admissible additions gives $G\in\I$. In both
cases involving torsion, homotopy invariance gives $\widetilde H_1(G;\ZZ)\cong
\ZZ/3$ and trivial reduced integral homology in every other degree.
Thus every integrally acyclic $K_4$-free graph of order $14$ lies in
$\I$. The reverse inclusion follows from~\eqref{eq:table-chain}.

The enumeration of all connected $K_4$-free graphs of order $14$ gives
\[
 \RowFourteenAcyclic-\RowFourteenIS=633
\]
isomorphism classes in $\mathcal A_{\mathbb F_2}\setminus\IS$. By
Theorem~\ref{thm:classification14}, eight of these classes belong to
$\I$, leaving $625$ with first integral homology $\ZZ/3$. The two enumerations use the same classifier. As a check on their
agreement, applying the neighbourhood criteria to all $633$
representatives in the full enumeration recovers the $51$ classes
found under the degree bounds. Computing integral homology directly
on the $625$ representatives gives $\ZZ/3$ in degree one and trivial
reduced homology in every other degree.
\end{proof}

\subsection{Minimum orders and dismantling levels}

The classification through order $13$ gives lower bounds without a
restriction on cliques. The graphs $\Gt$ and $\Ga$ attain these bounds,
both for graph reductions and for the associated flag complexes.

\begin{corollary}\label{cor:minimality}
The minimum order of a graph in $\mathcal A_{\mathbb F_2}\setminus\IS$
is $13$. The minimum order in each of $\I\setminus\IS$ and
$\mathcal A_{\ZZ}\setminus\IS$ is $14$. The minimum order of a graph
$G\in\IS$ with an admissible edge $e$ such that $G-e\notin\IS$
is also $14$. The minimum number of vertices of a contractible,
non-collapsible flag simplicial complex is $14$.
\end{corollary}
\begin{proof}
Theorem~\ref{thm:classification13} gives lower bounds of $13$ for
the order of a graph in $\mathcal A_{\mathbb F_2}\setminus\IS$ and
$14$ for the order of a graph in $\I\setminus\IS$. The graphs $\Gt$
and $\Ga$ attain the respective bounds. By
Proposition~\ref{prop:table-integral}, $14$ is also a lower bound
for the order of a graph in $\mathcal A_{\ZZ}\setminus\IS$, again
attained by $\Ga$. If $G\in\IS$ and an admissible edge $e$ satisfies
$G-e\notin\IS$, then $G-e\in\I\setminus\IS$ because admissible edge
deletions preserve I-contractibility. Thus $G-e$, and hence $G$, has
at least $14$ vertices. The one-edge extension of $\Ga$ in
Theorem~\ref{thm:B} attains this bound.

If a contractible flag complex $K=\cl(G)$ has at most $13$ vertices,
then $G\in\I$ by \cite[Corollary~1.5, p.~2]{DochtermannMatsushita2025}.
Theorem~\ref{thm:classification13} gives $G\in\IS$, so $K$ is
collapsible by \cite[Theorem~1]{EFH2022}. The contractible flag complex
$\cl(\Ga)$ has $14$ vertices and is non-collapsible by
Proposition~\ref{prop:noncollapsible}.
\end{proof}

For the $K_4$-free graphs of order $14$, the equality $\Dis{1}=\ISV$
follows from the neighbourhood criterion. At each step of a
vertex reduction witnessing membership in $\ISV$, the neighbourhood
of the deleted vertex is triangle-free and belongs to $\ISV$, and
hence is a nonempty tree. Such a tree is $0$-dismantlable, so the
same sequence gives a $1$-dismantling.

At order $12$ the equality $\Dis{1}=\Dis{2}$ first fails. By
subtraction in Table~\ref{tab:card}, $\Dis{2}\setminus\Dis{1}$ has
$20$ isomorphism classes at order $12$ and
$\NumTwoDisNotOneDisThirteen$ at order $13$. Table~2 of the supplementary material lists the
$20$ graphs of order $12$. The least order at which $\IS$ and $\ISV$
differ is $11$, with twelve classes in $\IS\setminus\ISV$
\cite[Theorem~11]{DEFH2023}. For orders $3$ through $11$, all $63$
counts for connected graphs, $\mathcal A_{\mathbb F_2}$, $\I$, $\IS$,
$\ISV$, $\Dis{1}$ and $\Dis{0}$ agree with \cite[Table~1]{DEFH2023}.

For the column counting $0$-dismantlable graphs there is also a
characterisation in terms of pursuit on a graph. Its entries count
precisely the cop-win graphs
\cite{NowakowskiWinkler1983}. At order $10$ we find $\num{2258314}$
such graphs up to isomorphism, in agreement with \cite[Table~1]{DEFH2023}
and one more than the $\num{2258313}$ reported by Baird et al.\
\cite[p.~72]{BairdEtAl2014}. Successive deletion of dominated vertices
and a direct test for a winning cop strategy select the same
$\num{2258314}$ representatives among the $\num{11716571}$ connected
graphs of this order.

Adjoining a universal vertex gives a bijection from the isomorphism
classes of graphs on $n-1$ vertices to those of cones on $n$ vertices.
Deleting any universal vertex gives the inverse, since any two
universal vertices are interchanged by an automorphism. This accounts
for the cone counts through order $13$ \cite{OEISA000088}. At these
orders the numbers of connected graphs also agree with
\cite{OEISA001349}. A cone is $K_4$-free precisely when deleting its
apex leaves a triangle-free graph. Thus the $\RowFourteenCone$ cones
in the last row correspond to the triangle-free graphs on $13$ vertices
\cite{OEISA006785}.

\appendix
\labelformat{section}{Appendix~#1}
\labelformat{subsection}{Appendix~#1}
\setcounter{table}{0}
\renewcommand{\thetable}{\thesection.\arabic{table}}
\renewcommand{\theHtable}{\theHsection.\arabic{table}}
\section{The disc identification and the explicit graph reductions}\label{app:certificates}

The explicit descriptions below use the vertex labels of
Table~\ref{tab:witness-adjacency}. For $\Gt$ they give the disc
whose boundary determines the attaching map. For the graphs on
fourteen vertices they give both the obstruction to a first deletion
and the reductions that become possible after one edge is added.

\subsection{A triangulated disc for the thirteen-vertex graph}
\label{app:disc-certificates}

In Figure~\ref{fig:G13} each occurrence of a label on the boundary
of the disc represents a separate vertex before the identifications.
To distinguish them, write $b_0,\ldots,b_{11}$ for the boundary
vertices in clockwise order, with $b_0$ at the top. The remaining
vertices are $0,1,\ldots,7,12$. Let $D$ be the simplicial complex
whose triangles are listed in Table~\ref{tab:G13-disc-triangles}.
The map $q:D\to\cl(\Gt)$ fixes the labels of the interior vertices
and is given on the boundary by
\[
 b_{4r}\longmapsto10,\qquad
 b_{4r+1}\longmapsto9,\qquad
 b_{4r+2}\longmapsto11,\qquad
 b_{4r+3}\longmapsto8\qquad(0\le r\le2).
\]

\begin{table}[htbp]
\centering\small
\caption{Triangles of the disc $D$ in the construction of $\cl(\Gt)$.}\label{tab:G13-disc-triangles}
\begin{tabular}{@{}rlrlrlrl@{}}
\toprule
$j$ & $T_j$ & $j$ & $T_j$ & $j$ & $T_j$ & $j$ & $T_j$ \\
\midrule
1 & $\{0,4,b_{4}\}$ & 8 & $\{1,4,12\}$ & 15 & $\{2,7,b_{9}\}$ & 22 & $\{3,7,12\}$ \\
2 & $\{0,4,12\}$ & 9 & $\{1,6,b_{1}\}$ & 16 & $\{2,7,12\}$ & 23 & $\{3,b_{11},b_{0}\}$ \\
3 & $\{0,5,b_{6}\}$ & 10 & $\{1,6,12\}$ & 17 & $\{2,b_{7},b_{8}\}$ & 24 & $\{3,b_{11},b_{10}\}$ \\
4 & $\{0,5,12\}$ & 11 & $\{1,b_{3},b_{2}\}$ & 18 & $\{2,b_{9},b_{8}\}$ & 25 & $\{4,b_{3},b_{4}\}$ \\
5 & $\{0,b_{5},b_{4}\}$ & 12 & $\{1,b_{1},b_{2}\}$ & 19 & $\{3,6,b_{0}\}$ & 26 & $\{5,b_{7},b_{6}\}$ \\
6 & $\{0,b_{5},b_{6}\}$ & 13 & $\{2,5,b_{7}\}$ & 20 & $\{3,6,12\}$ & 27 & $\{6,b_{1},b_{0}\}$ \\
7 & $\{1,4,b_{3}\}$ & 14 & $\{2,5,12\}$ & 21 & $\{3,7,b_{10}\}$ & 28 & $\{7,b_{9},b_{10}\}$ \\
\bottomrule
\end{tabular}

\end{table}

Begin with $T_1$ and attach the other triangles in the following order
\[
\begin{gathered}
2,4,3,5,6,8,7,10,9,11,12,14,13,16,15,17,18,\\
20,19,22,21,23,24,25,26,27,28.
\end{gathered}
\]
At each attachment in this order, Table~\ref{tab:G13-disc-triangles}
gives the intersection with the preceding union. It consists of a
single boundary edge or two consecutive boundary edges, together with
their endpoints. Gluing the new triangle along this arc preserves the
disc, so induction gives a triangulated disc $D$ with boundary cycle
$b_0,b_1,\ldots,b_{11},b_0$.

We now check that the boundary identifications account for every
identification made by $q$. Tables~\ref{tab:G13-disc-triangles} and~\ref{tab:witness-adjacency}
show that the triangles of $D$ have distinct images and cover all $28$ triangles
of $\cl(\Gt)$. Their images contain every edge of $\Gt$ and,
since $\Gt$ is $K_4$-free, they cover its entire clique complex.
The interior vertex labels are distinct and lie outside
$\{8,9,10,11\}$. Each interior edge likewise has a distinct image
with an endpoint outside the boundary four-cycle. Consequently,
only boundary vertices with the same label and the corresponding
boundary edges are identified. This proves the quotient description
in Proposition~\ref{prop:disc}. The boundary cycle traverses the
four-cycle three times with the same orientation, so the attaching
map has degree three.

\subsection{Obstructions to admissible deletions}

To decide whether a deletion can begin a strong reduction, we need
the homology of its neighbourhood in the original graph.
For each graph $G$ in Table~\ref{tab:witness-links}, the neighbourhood
$\Nb_G(v)$ of any vertex $v$ is connected and triangle-free, so its
first Betti number is $|E(\Nb_G(v))|-|V(\Nb_G(v))|+1$. The middle
column gives these numbers in increasing order of the vertex labels.
For an edge $uv$, each vertex of $\Nb_G(u,v)$ is a common neighbour
of its endpoints. Since $G$ is $K_4$-free, no two of these vertices
are adjacent, so the order of $\Nb_G(u,v)$ determines its number
of components. The last column lists all orders that occur.

Every entry in the middle column is positive, and every order in the
last column is at least two. Thus each vertex neighbourhood has
nonzero first homology and each edge neighbourhood has nonzero
reduced zeroth homology. Thus no vertex or edge has the neighbourhood in $\IS$ required
for a first admissible deletion.

\begin{table}[htbp]
\centering\small
\caption{Neighbourhoods of vertices and common neighbourhoods of edge endpoints.
The first Betti numbers are independent of the coefficient field.}\label{tab:witness-links}
\begin{tabular}{@{}lll@{}}
\toprule
Graph & First Betti numbers of vertex neighbourhoods & Common neighbourhood orders\\
\midrule
$\Gt$ & $1,1,1,1,1,1,1,1,2,2,2,2,1$ & $2,3$\\
$\Ga$ & $1,1,1,1,1,2,1,1,1,2,2,1,2,2$ & $2,3$\\
$\Gb$ & $1,1,1,1,1,2,1,2,2,2,1,1,1,3$ & $2,3,4$\\
$\Gc$ & $2,1,1,1,1,1,1,2,1,2,2,2,1,2$ & $2,3$\\
\bottomrule
\end{tabular}
\end{table}

\subsection{Reductions of the one-edge extensions}

After the indicated edge is added to each of $\Ga,\Gb,\Gc$, a
strong reduction is possible. Tables~\ref{tab:Ga-reduction}
and~\ref{tab:variant-reductions} give the deletions in order and the
neighbourhood just before each deletion. A single vertex label denotes
a one-vertex graph, and a sequence joined by dashes denotes the path
through those vertices in the given order. At step $13$ for
$\Gb+\{0,3\}$, two paths meet at vertex $8$, so their union is
again a tree. In particular, every listed neighbourhood is a
nonempty tree and belongs to $\IS$. The deletions leave vertex $13$
in each of the three extensions.

\begin{table}[htbp]
\centering\small
\caption{The complete reduction of $\Ga+\{0,1\}$ used in
Theorem~\ref{thm:B}.}\label{tab:Ga-reduction}
\begin{tabular}{@{}rll@{}}
\toprule
Step & Delete & Current neighbourhood \\
\midrule
1 & $e(0,5)$ & $10\text{--}1\text{--}12$ \\
2 & $e(0,10)$ & $1\text{--}13$ \\
3 & $e(0,12)$ & $1\text{--}8$ \\
4 & $e(8,12)$ & $1$ \\
5 & $e(1,12)$ & $5$ \\
6 & $v(5)$ & $1\text{--}10\text{--}11\text{--}9\text{--}12$ \\
7 & $v(10)$ & $1\text{--}13\text{--}3\text{--}6\text{--}11$ \\
8 & $v(3)$ & $6\text{--}12\text{--}9\text{--}13$ \\
9 & $v(6)$ & $11\text{--}7\text{--}2\text{--}12$ \\
10 & $v(2)$ & $7\text{--}13\text{--}9\text{--}12$ \\
11 & $v(7)$ & $11\text{--}8\text{--}1\text{--}13$ \\
12 & $v(1)$ & $8\text{--}0\text{--}13$ \\
13 & $v(0)$ & $8\text{--}4\text{--}13$ \\
14 & $v(4)$ & $8\text{--}11\text{--}9\text{--}13$ \\
15 & $v(8)$ & $11$ \\
16 & $v(11)$ & $9$ \\
17 & $v(12)$ & $9$ \\
18 & $v(9)$ & $13$ \\
\bottomrule
\end{tabular}

\end{table}

\begin{table}[htbp]
\centering\small
\caption{Reductions of $\Gb+\{0,3\}$ and $\Gc+\{2,3\}$.
A dash marks the end of the shorter reduction.}\label{tab:variant-reductions}
\begin{tabular}{@{}rllll@{}}
\toprule
& \multicolumn{2}{c}{$\Gb+\{0,3\}$} & \multicolumn{2}{c}{$\Gc+\{2,3\}$} \\
Step & Delete & Current neighbourhood & Delete & Current neighbourhood \\
\midrule
1 & $e(0,9)$ & $12\text{--}3\text{--}13$ & $e(2,7)$ & $10\text{--}3\text{--}13$ \\
2 & $e(0,13)$ & $3\text{--}7$ & $e(2,10)$ & $3\text{--}11$ \\
3 & $e(9,12)$ & $3$ & $e(2,13)$ & $3\text{--}6$ \\
4 & $e(3,9)$ & $13$ & $e(6,13)$ & $3$ \\
5 & $e(3,13)$ & $7$ & $e(3,13)$ & $7$ \\
6 & $e(7,13)$ & $2$ & $v(7)$ & $3\text{--}10\text{--}12\text{--}0\text{--}13$ \\
7 & $v(2)$ & $7\text{--}11\text{--}6\text{--}12\text{--}8\text{--}13$ & $v(10)$ & $3\text{--}11\text{--}5\text{--}1\text{--}12$ \\
8 & $v(6)$ & $11\text{--}4\text{--}10\text{--}3\text{--}12$ & $v(1)$ & $5\text{--}13\text{--}4\text{--}8\text{--}12$ \\
9 & $v(3)$ & $10\text{--}7\text{--}0\text{--}12$ & $v(4)$ & $8\text{--}11\text{--}9\text{--}0\text{--}13$ \\
10 & $v(0)$ & $7\text{--}5\text{--}12$ & $v(0)$ & $12\text{--}9\text{--}5\text{--}13$ \\
11 & $v(4)$ & $10\text{--}8\text{--}13\text{--}9\text{--}11$ & $v(8)$ & $11\text{--}2\text{--}6\text{--}12$ \\
12 & $v(7)$ & $10\text{--}5\text{--}11$ & $v(2)$ & $6\text{--}3\text{--}11$ \\
13 & $v(5)$ & $11\text{--}1\text{--}8\text{--}10,\ 8\text{--}12$ & $v(3)$ & $6\text{--}9\text{--}11$ \\
14 & $v(1)$ & $8\text{--}13\text{--}9\text{--}11$ & $v(6)$ & $9\text{--}12$ \\
15 & $v(10)$ & $8$ & $v(11)$ & $5\text{--}9$ \\
16 & $v(11)$ & $9$ & $v(12)$ & $9$ \\
17 & $v(9)$ & $13$ & $v(9)$ & $5$ \\
18 & $v(12)$ & $8$ & $v(5)$ & $13$ \\
19 & $v(8)$ & $13$ & --- & --- \\
\bottomrule
\end{tabular}

\end{table}

\FloatBarrier
\section*{Data and software availability}
The accompanying supplement, \emph{Finite closures and graph certificates}, develops the admissible edge additions of
Section~\ref{sec:fourteen} and describes the surfaces obtained by
unfolding the associated complexes. The enumeration programs and
records, graph catalogues, reduction sequences and homology
calculations are available in release \texttt{v1.0.0} of the
companion repository~\cite{GraphData2026} at
\href{https://github.com/Geometric-and-Combinatorial-Structures/graph-reductions-flag-complexes}{\nolinkurl{github.com/Geometric-and-Combinatorial-Structures/}\allowbreak\nolinkurl{graph-reductions-flag-complexes}}. The supplementary archive \texttt{cop-win-10.zip}
contains two programs for counting cop-win graphs on ten vertices,
one by vertex domination and the other by optimal play.

\section*{Acknowledgements}
The authors thank the \'Area de C\'omputo de Alto Rendimiento de la
Universidad de Sonora (ACARUS) for computing time on the Yuca cluster
and for technical support.

\section*{Declarations}
\noindent\textbf{Funding.}
This work was supported by the Universidad de Sonora through the internal
research funding call of the Direcci\'on de Apoyo a Docentes,
Investigaci\'on y Posgrado (DADIP), project USO315010378. The University's institutional role was limited to funding and the technical support provided by ACARUS.

\noindent\textbf{Author contributions.}
The authors are listed alphabetically and contributed equally to this work.

\noindent\textbf{Competing interests.}
The authors declare that they have no known competing financial interests
or personal relationships that could have appeared to influence the work
reported in this paper.

\section*{Declaration of generative AI and AI-assisted technologies in the manuscript preparation process}
The authors used Codex (OpenAI) and Claude (Anthropic) to draft and
revise the exposition, arguments, tables and software, to check sources
and to examine the arguments for possible errors. The authors retain
responsibility for the mathematical content and the final manuscript.

\bibliographystyle{elsarticle-num-sorted}
\bibliography{refs}
\end{document}